\documentclass[11pt]{amsart}
\usepackage[T1]{fontenc}
\usepackage{lmodern}
\usepackage[margin=1.05in]{geometry}
\usepackage{amsmath,amssymb,amsthm,mathtools}
\usepackage{microtype}
\usepackage[colorlinks=true,linkcolor=blue,citecolor=blue,urlcolor=blue]{hyperref}
\newtheorem{theorem}{Theorem}[section]
\newtheorem{proposition}[theorem]{Proposition}
\newtheorem{lemma}[theorem]{Lemma}
\newtheorem{corollary}[theorem]{Corollary}
\theoremstyle{definition}

\newtheorem{example}[theorem]{Example}
\theoremstyle{remark}
\newtheorem{remark}[theorem]{Remark}
\newcommand{\R}{\mathbb R}
\newcommand{\G}{G(n,k)}

\newcommand{\nuu}{\,d\nu(E)}
\newcommand{\supp}{\operatorname{supp}}
\newcommand{\tr}{\operatorname{tr}}

\newcommand{\BL}{\mathcal B_\nu}
\newcommand{\RB}{\mathcal R_\nu}
\newcommand{\RBone}{\mathcal R_\nu^{\,1}}
\newcommand{\Mass}{\mathcal M_\nu}
\newcommand{\cT}{\mathcal T_S}

\DeclareMathOperator{\Cov}{Cov}

\numberwithin{equation}{section}
\title[Continuous Grassmannian Brascamp--Lieb inequalities]{Continuous Grassmannian Brascamp--Lieb inequalities and their equality cases}
\author[A.-J. Li]{Ai-Jun Li}
\address[A.-J. Li]{School of Science,
	Zhejiang University of Science and Technology, Hangzhou 310023, China}
\email{liaijun72@163.com}	
\author[Q. Huang]{Qingzhong Huang}
\address[Q. Huang]{College of data science, Jiaxing University, Jiaxing 314001,	China} \email{hqz376560571@163.com} 
\date{}
\begin{document}
\begin{abstract}
In this paper, we establish the forward and reverse geometric Brascamp--Lieb inequalities for isotropic measures on a Grassmann manifold, and characterize their equality conditions. The proofs use transportation and the Grassmannian Ball--Barthe inequality of Li, Xi, and Zhang. The functions form a jointly continuous family whose fiber integrals are positive and continuous. The equality conditions are described by an orthogonal decomposition determined by the supporting projections. On the common eigenspaces, the forward inequality permits arbitrary independent factors, whereas the reverse inequality requires log-concavity for factors shared by a nontrivial parameter distribution. On the orthogonal complement, the extremizers have Gaussian factors. Their centers are compatible projections in the forward inequality and may vary continuously between fibers in the reverse inequality. Gaussian convolution preserves extremizers, and the original factors are recovered at a fixed positive time using characteristic functions. As applications, we obtain Gaussian rigidity results and volume bounds for bodies defined by projection integrals and barycentric infima. These bounds recover the reverse isoperimetric inequalities of Li, Xi, and Zhang for bodies associated with Grassmannian cosine transforms, together with their equality conditions away from the endpoint cases. We also establish a forward extension under weaker mass assumptions and give a mass-escape example for the reverse inequality.
\end{abstract}
\thanks{The first author was supported by NSFC (Grant No. 12571146 and 12231006). The second author was supported by  NSFC (Grant No. 11701219) and the Qin Shen Scholar Program of Jiaxing University.}
\maketitle
\section{Introduction}

The general Brascamp--Lieb inequalities provide a common framework for several fundamental integral inequalities. In particular, suitable choices of the  linear maps and exponents yield H\"{o}lder's inequality, the Loomis--Whitney inequality, and, by duality, Young's convolution inequality. The sharp constants in these inequalities reflect the interaction between the linear maps and the geometry of their underlying spaces. Brascamp and Lieb \cite{BL} established the sharp form of Young's inequality and its multilinear generalizations. Lieb \cite{Lieb} subsequently showed that the optimal constants in the general Brascamp--Lieb inequality are determined by Gaussian functions. The finiteness of these constants and the structure of extremizing data were studied systematically by Bennett, Carbery, Christ, and Tao \cite{BCCT}. These results show the role of Gaussian functions in the determination of the optimal constants. They also show that the existence and structure of extremizers depend on the underlying linear maps and exponents.

A particularly useful form in convex geometry is the geometric Brascamp--Lieb inequality. Let $E_1,\ldots,E_N$ be linear subspaces of $\R^n$, let $P_{E_i}$ denote orthogonal projection onto $E_i$, and suppose that positive numbers $c_1,\ldots,c_N$ satisfy
\begin{equation}\label{eq:finite-isotropy}
	\sum_{i=1}^N c_iP_{E_i}=I_n.
\end{equation}
Then, for nonnegative integrable functions $f_i$ on $E_i$,
\begin{equation}\label{eq:finite-forward}
	\int_{\R^n}\prod_{i=1}^N f_i(P_{E_i}x)^{c_i}\,dx
	\le \prod_{i=1}^N\left(\int_{E_i}f_i\right)^{c_i}.
\end{equation}
Under this identity decomposition, the optimal constant in the inequality is one. Ball's work \cite{Ball} demonstrated the geometric significance of this normalization through volume-ratio estimates and a reverse isoperimetric inequality. Barthe \cite{Barthe} proved the corresponding reverse inequality,
\begin{equation}\label{eq:finite-reverse}
	\int_{\R^n}^{*}\sup_{\substack{y_i\in E_i\\\sum_i c_i y_i=x}}
	\prod_{i=1}^N f_i(y_i)^{c_i}\,dx
	\ge \prod_{i=1}^N\left(\int_{E_i}f_i\right)^{c_i},
\end{equation}
where $\int^{*}$ denotes the outer integral. His transportation argument connects both inequalities with a determinant inequality arising from the identity decomposition. Thus the product in \eqref{eq:finite-forward} and the supremum in \eqref{eq:finite-reverse} provide two ways of combining functions defined on the subspaces.

The equality conditions for \eqref{eq:finite-forward} and \eqref{eq:finite-reverse} depend on the geometry of the subspaces. An extremizer may contain non-Gaussian factors associated with common eigenspaces of the supporting projections, and the conditions on these factors differ in the two inequalities. Valdimarsson \cite{Valdimarsson} characterized the extremizers of \eqref{eq:finite-forward}. The equality conditions for the geometric reverse inequality were obtained by B\"or\"oczky, Kalantzopoulos, and Xi \cite{BKX}. In the reverse inequality, a shared independent factor must be log-concave, while the Gaussian factors may be translated separately in the fibers. In the forward inequality, the Gaussian centers must be projections of a common vector. A broader finite-dimensional framework combining forward and reverse Brascamp--Lieb inequalities was developed by Courtade and Liu \cite{CL}. The present paper concerns the geometric pair \eqref{eq:finite-forward}--\eqref{eq:finite-reverse} with a continuously indexed family of subspaces.

A continuous version replaces the finite sum in \eqref{eq:finite-isotropy} by an integral with respect to an isotropic measure. For measures on lines, Barthe \cite{BartheContinuous} proved both the forward and reverse inequalities. Isotropic measures also appear in volume inequalities. Lutwak, Yang, and Zhang \cite{L1} extended Ball and Barthe's reverse isoperimetric inequalities to unit balls of subspaces of $L_p$, using even isotropic measures. They also established sharp inequalities and their equality conditions for arbitrary isotropic measures \cite{LYZ}. These results lead to the question of what happens when lines are replaced by higher-dimensional subspaces. A related direction, in which the geometric normalization is allowed to be approximate, was considered by Brazitikos and Giannopoulos \cite{BG}, who established continuous versions of approximate geometric Brascamp--Lieb inequalities.

Let $G(n,k)$ be the set of $k$-dimensional linear subspaces of $\R^n$. For a positive Borel measure $\nu$ on $G(n,k)$, the isotropy condition is
\[
 \int_{G(n,k)}P_E\,d\nu(E)=I_n.
\]
Li, Xi, and Zhang \cite{LXX} developed volume inequalities for convex bodies arising from cosine transforms on Grassmann manifolds. A principal ingredient of their approach is a Grassmannian Ball--Barthe determinant inequality. Li, Cao, and Huang \cite{LCH} subsequently studied a matrix form of this inequality. The Grassmannian Ball--Barthe inequality compares the determinant of an average of positive operators on the subspaces with the exponential of the integral of their logarithmic determinants. In a transportation argument, the operators on the subspaces are derivatives of the transport maps, while their average is the derivative of the induced map on the ambient space. The equality condition then relates the fiber derivatives to a common ambient operator. This connection is the starting point of our proofs.

Continuous Brascamp--Lieb inequalities also appear in the stochastic framework of Cordero-Erausquin and Maurey \cite{CEM}. In particular, their Section~5 treats generalized geometric forward and reverse inequalities through operators between spaces of square-integrable functions. Their work provides a broader context for passing from finite index sets to measure spaces. More recently, Li \cite{Li2026} studied the continuous Grassmannian forward inequality by measure approximation and heat flow. The measure-approximation argument preserves the isotropy constraint while passing from finite data to a measure on the Grassmannian. These works provide the continuous inequalities and methods on which the present study builds. Our emphasis is on a direct Grassmann Ball--Barthe transportation treatment of the two inequalities, together with their equality conditions on a common class of continuous families.

Heat flow provides another approach to these inequalities. Carlen, Lieb, and Loss \cite{CLL} used a nonlinear heat flow to study sharp integral inequalities and their optimizers. Bennett, Carbery, Christ, and Tao \cite{BCCT} developed heat-flow arguments for general Brascamp--Lieb data, and Barthe and Huet \cite{BH} gave parallel semigroup proofs of the geometric forward and reverse inequalities. In the present paper, transportation proves the inequalities, whereas Gaussian convolution is used to analyze equality. This division is useful because the inequality remains valid for smooth approximations, whereas equality need not be preserved.

In this paper, we establish the forward and reverse Grassmannian inequalities for jointly continuous families whose positive finite fiber masses depend continuously on the subspace. Equivalently, after normalization, we work with continuous families of probability densities on the tautological bundle over the support of $\nu$. We characterize equality in terms of the complete common $0/1$ eigenspaces of the supporting projections and their orthogonal complement. On the latter space, the extremizers have a common Gaussian covariance commuting with the projections. On the common eigenspaces, the forward inequality allows arbitrary continuous probability factors, while the reverse inequality requires log-concavity precisely when a factor is shared by a nontrivial parameter distribution. The Gaussian centers in the forward case are projections of one ambient vector. In the reverse case, a continuous translation may be chosen separately in each fiber. For finitely supported measures, these descriptions agree with the classifications in \cite{Valdimarsson,BKX}.

In the main theorems, the functions form a jointly continuous family whose fiber integrals are positive and continuous. In Appendix~\ref{sec:extension}, we extend the forward inequality and its equality characterization to a broader class of families. Appendix~\ref{sec:escape} gives an example illustrating why the same extension cannot be made for the reverse inequality without additional assumptions.

As applications, we first show that, for $1\le k<n$, if the supporting projections have no nontrivial common invariant subspace, equality can hold only for Gaussian functions. We also characterize equality when the functions have a common radial profile. Moreover, we obtain volume bounds for convex bodies defined by projection integrals and the corresponding barycentric infima. For $1<p<\infty$, the two families of bodies with conjugate exponents are related by polarity. For finite $p\ge1$ with $p\ne2$, equality in either volume bound is characterized by an orthogonal configuration of subspaces. These bounds recover the reverse isoperimetric inequalities for the bodies associated with Grassmannian cosine transforms established by Li, Xi, and Zhang \cite{LXX}; see Corollary~\ref{cor:lxx-comparison}. The finite-exponent equality conditions follow from our functional equality conditions. The endpoint inequalities are obtained by taking limits, while the remaining endpoint equality conditions are stated with reference to \cite{LXX}.

\subsection{Notation and admissible families}
We work in the Euclidean space $\R^n$. For $1\le k\le n$,  let $\G$ be the Grassmann manifold of $k$-dimensional subspaces of $\R^n$, and let $\nu$ be a finite positive Borel measure on $\G$. Write $S=\supp\nu$. For $E\in\G$, let $P_E:\R^n\to\R^n$ be the orthogonal projection onto $E$. All integrals on $E$ use the induced Lebesgue measure. We assume that
\begin{equation}\label{eq:isotropy}
	\int_S P_E\nuu=I_n.
\end{equation}
Taking traces in \eqref{eq:isotropy} gives $\nu(S)=n/k$. Thus, if $\mu$ is an isotropic probability measure with $\int P_E\,d\mu=(k/n)I_n$, then $\nu=(n/k)\mu$ satisfies \eqref{eq:isotropy}.

Write $\cT=\{(E,y):E\in S,\ y\in E\}$ for the tautological bundle restricted to $S$. A family $f=(f_E)$ is called \emph{admissible} if the function $f:\cT\to[0,\infty)$ is jointly continuous and
\begin{equation}\label{eq:massclass}
	m_E=\int_Ef_E(y)\,dy\in(0,\infty),\qquad E\mapsto m_E
	\text{ is continuous on }S.
\end{equation}
Put $q_E=f_E/m_E$. Then $\int_Eq_E=1$ for every $E\in S$. By Scheff\'e's lemma in local orthonormal frames and the compactness of $S$, the densities are uniformly tight:
\begin{equation}\label{eq:tightness}
	\lim_{R\to\infty}\sup_{E\in S}\int_{\{|y|>R\}\cap E}q_E(y)\,dy=0.
\end{equation}
Set
\[
\begin{gathered}
	F_f(x)=\exp\int_S\log f_E(P_Ex)\nuu,\\
	\BL(f)=\int_{\R^n}F_f(x)\,dx,\qquad
	\Mass(f)=\exp\int_S\log m_E\nuu.
\end{gathered}
\]
We use $\log0=-\infty$ and $e^{-\infty}=0$. For each fixed $x$, the positive part of
$\log f_E(P_Ex)$ is bounded, so $F_f$ is well defined.

A section is a $\nu$-measurable field $v(E)\in E$, considered up to equality almost everywhere. Define the reverse envelope by
\begin{equation}\label{eq:reverse-envelope}
	\RB(f)(x)=
	\sup_{\substack{v\in L^\infty(\nu;\R^n),\ v(E)\in E\\
			\int_Sv(E)\nuu=x}}
	\exp\int_S\log f_E(v(E))\nuu.
\end{equation}
The constraint is nonempty by \eqref{eq:isotropy}, since one may take
$v(E)=P_Ex$. We also write $\RBone(f)$ for the corresponding envelope
over Bochner-integrable sections satisfying
$\int_S\log^+f_E(v(E))\,d\nu(E)<\infty$. More explicitly,
\[
 \RBone(f)(x)=
 \sup_{\substack{v\in L^1(\nu;\R^n),\ v(E)\in E\ \nu\text{-a.e.}\\
                 \int_Sv(E)\,d\nu(E)=x,\ 
                 \int_S\log^+f_E(v(E))\,d\nu(E)<\infty}}
 \exp\int_S\log f_E(v(E))\,d\nu(E).
\]
Both envelopes are universally measurable
(see Proposition~\ref{prop:measurability}) and agree with the usual
finite-support expression. The same inequality and equality
characterization hold for both envelopes.

\subsection{The subspaces determining equality}
The equality conditions depend on the geometry of the support $S$. To describe them, for each $\varepsilon:S\to\{0,1\}$, define
\[
F_\varepsilon=\{z:P_Ez=\varepsilon(E)z\text{ for every }E\in S\}.
\]
Equivalently,
\[
F_\varepsilon
=
\left(\bigcap_{\varepsilon(E)=1}E\right)
\cap
\left(\bigcap_{\varepsilon(E)=0}E^\perp\right).
\]
An empty intersection is understood as $\R^n$. The distinct nonzero spaces $F_\varepsilon$ are pairwise orthogonal. List them as $F_1,\ldots,F_s$, where $0\le s\le n$, without further subdivision, and set
\begin{equation}\label{eq:decomposition}
	D=(F_1\oplus\cdots\oplus F_s)^\perp,\qquad E_D=E\cap D.
\end{equation}
For every $E\in S$, we have
\begin{equation}\label{eq:block-identities}
	\begin{gathered}
		E=E_D\oplus\bigoplus_{F_j\subset E}F_j,\qquad
		S_j=\{E\in S:F_j\subset E\},\\
		\nu(S_j)=1,\qquad \int_SP_{E_D}\nuu=I_D.
	\end{gathered}
\end{equation}
Here $P_{E_D}$ denotes the orthogonal projection onto $E_D$; in the last identity it is restricted to $D$. To see the identities, take a unit vector $u\in F_j$. For $E\in S$,
\[
 P_Eu=\begin{cases}u,&E\in S_j,\\0,&E\notin S_j.\end{cases}
\]
Thus isotropy gives
\[
 \nu(S_j)=\int_S\langle P_Eu,u\rangle\,d\nu(E)=|u|^2=1.
\]
Each $P_E$ preserves the spaces $F_j$ and $D$. Hence
\[
 P_E|_D=P_{E_D}|_D,\qquad
 \int_S P_{E_D}|_D\,d\nu(E)=I_D.
\]
The function $E\mapsto|P_Eu|^2$ is continuous on $S$ and takes only the values $0$ and $1$. Therefore $S_j$ is both open and closed relative to $S$.

We call $F_j$ \emph{shared} if the probability measure $\nu|_{S_j}$ is not a point mass. Equivalently, there is a Borel set $A\subset S_j$ such that $0<\nu(A)<1$. The space $D$ may contain nontrivial common reducing subspaces of the projections $P_E$, $E\in S$.

\begin{theorem}[Forward inequality and equality]\label{thm:forward}
	For an admissible family $f$,
	\begin{equation}\label{eq:forward-main}
		\BL(f)\le\Mass(f).
	\end{equation}
	The constant $1$ is optimal. Equality holds if and only if there exist a positive definite symmetric operator $A:D\to D$ commuting with every $P_{E_D}$, a vector $a\in D$, continuous probability densities $h_j$ on $F_j$, and positive continuous numbers $\theta_E$, such that
	\begin{equation}\label{eq:forward-form}
		f_E(y)=\theta_E
		e^{-\pi\langle A(P_Dy-P_{E_D}a),\,P_Dy-P_{E_D}a\rangle}
		\prod_{F_j\subset E}h_j(P_{F_j}y)
		\quad(E\in S,\ y\in E).
	\end{equation}
	 Empty products and zero-dimensional Gaussian factors are $1$.
\end{theorem}

\begin{theorem}[Reverse inequality and equality]\label{thm:reverse}
	For an admissible family $f$,
	\begin{equation}\label{eq:reverse-main}
		\int_{\R^n}\RB(f)(x)\,dx\ge\Mass(f).
	\end{equation}
	The constant $1$ is optimal. Equality holds if and only if there exist a positive definite symmetric operator $A:D\to D$ commuting with every $P_{E_D}$, a continuous section $b_E\in E$, continuous probability densities $h_j$ on $F_j$, and positive continuous numbers $\theta_E$, such that
	\begin{equation}\label{eq:reverse-form}
		f_E(y)=\theta_E
		e^{-\pi\langle A(P_Dy-P_Db_E),\,P_Dy-P_Db_E\rangle}
		\prod_{F_j\subset E}h_j\big(P_{F_j}(y-b_E)\big)
		\quad(E\in S,\ y\in E),
	\end{equation}
	where $h_j$ is log-concave whenever $F_j$ is shared.
	The same assertion, including the equality characterization, holds for $\RBone$.
\end{theorem}


The forward formula has one Gaussian center whose projections occur in all fibers. The reverse
formula permits a separate continuous translation in each fiber. An independent factor occurring at a single unit atom need not be log-concave, since the inequality on that block reduces to Fubini's identity.

When $\nu=\sum_{i=1}^N c_i\delta_{E_i}$, the two theorems reduce to the finite geometric inequalities with their equality conditions. The assumptions on the fiber masses are then automatic for positive integrable continuous functions. When \(k=n\), the Grassmann manifold \(G(n,n)\) consists of the single subspace \(\mathbb R^n\). Isotropy then gives \(\nu=\delta_{\mathbb R^n}\), and both inequalities reduce to identities. For $k<n$ and full support, the common Gaussian covariance is scalar and there are no independent factors; see Corollary~\ref{cor:full}.

\subsection{Outline of the proofs}
The Grassmannian Ball--Barthe inequality plays two roles in the proofs. Its determinant bound gives the Jacobian estimate, while its equality condition yields commutation with the supporting projections. A corresponding rigidity lemma shows that, if the Hessian of a sufficiently smooth function commutes with all the supporting projections, then the function is quadratic in the \(D\)-variable and can be written as a sum of functions depending separately on the variables in \(F_1,\ldots,F_s\). We establish these facts in Section~\ref{sec:structure}.

In Section~\ref{sec:transport},
we first establish the measurability of the reverse envelopes and the
continuous dependence of the transports between fixed balls. To prove
the reverse inequality, we use weak transports with bounded targets
and an area inequality for the gradient of their averaged potential. For each fixed positive regularization time, the transport maps and their derivatives depend continuously on \(E\) in local coordinates, uniformly on compact sets.

In Section~\ref{sec:inequalities}, we prove the forward and reverse
inequalities. We then study their equality conditions in
Sections~\ref{sec:forward-equality} and~\ref{sec:reverse-equality}.
Using Gaussian convolution, we show that an extremizer remains an
extremizer after regularization. This allows us to first characterize
the equality cases for smooth functions.
To recover the original functions, we fix a positive regularization
time. The factors in the resulting equality conditions are identified
with marginals of the original probability densities. We then use
characteristic functions to recover their product structure and show
that the factors agree across the fibers.

Finally, Section~\ref{sec:consequences} gives the rigidity and volume consequences. Appendix~\ref{sec:extension} proves the broader forward theorem, and Appendix~\ref{sec:escape} explains the restriction on the family in the reverse theorem.

\section{The determinant inequality and the common subspaces}\label{sec:structure}

We begin with the Grassmannian Ball--Barthe inequality and an operator formulation of its equality condition. We then establish a rigidity statement for functions whose Hessians commute with the supporting projections. Finally, we give a deconvolution lemma for recovering probability factors after Gaussian regularization.

For each $E\in S$, write $\iota_E:E\to\R^n$ for the natural inclusion.

\begin{lemma}[Grassmannian Ball--Barthe inequality]\label{lem:BB}
	For a continuous positive definite symmetric field $B_E:E\to E$, put $M=\int_S\iota_EB_EP_E\nuu$. Then
	\begin{equation}\label{eq:BB}
		\det M\ge\exp\int_S\log\det_E B_E\nuu.
	\end{equation}
	Equality holds if and only if
	\begin{equation}\label{eq:BB-equality}
		MP_E=P_EM,\qquad B_E=M|_E\quad(E\in S).
	\end{equation}
	The inequality also holds for measurable positive definite fields with
	$\int\tr_EB_E\,d\nu<\infty$ and
	$\int|\log\det_EB_E|\,d\nu<\infty$.
\end{lemma}
\begin{proof}
	The continuous inequality is Lemma 3.4 of \cite{LXX}. The original inequality is stated for a continuous positive definite
	field on $G(n,k)$. To apply it here, we first regard
	$\log B_E$ as a symmetric operator on $\R^n$ by setting
	\[
	H(E)=\iota_E(\log B_E)P_E,\qquad E\in S.
	\]
	Since $S$ is closed, $H$ admits a continuous symmetric extension
	$\widetilde H$ to $G(n,k)$. For $F\in G(n,k)$, define
	\[
	\widetilde B_F
	=\exp_F\bigl(P_F\widetilde H(F)|_F\bigr).
	\]
	Then $\widetilde B_F$ is positive definite, and
	$\widetilde B_E=B_E$ for every $E\in S$. The desired inequality
	follows from the original one, since $\nu$ is supported on $S$.

	To examine the equality condition, consider the functional
	\[
	\mathcal D(B)
	=
	\log\det M-\int_S\log\det_E B_E\,d\nu(E).
	\]
	By the Grassmannian Ball--Barthe inequality, $\mathcal D(B)\ge0$.
	Suppose that equality holds. Then $\mathcal D(B)=0$.
	Let $H_E:E\to E$ be a continuous family of symmetric operators, and
	consider the variation
	\[
	B_E(t)=B_E+tH_E.
	\]
	Since $B_E$ is positive definite and depends continuously on $E$, the
	compactness of $S$ gives a constant $c>0$ such that
	\[
	B_E\ge cI_E
	\qquad(E\in S).
	\]
	The family $H_E$ is uniformly bounded on $S$. Hence, for sufficiently
	small $|t|$, the operators $B_E(t)$ are positive definite for every
	$E\in S$. 
	Put
	\[
	M(t)=\int_S\iota_EB_E(t)P_E\,d\nu(E).
	\]
	Since $\mathcal D(B(t))\ge0$ and $\mathcal D(B(0))=0$, we have
	\[
	\left.\frac{d}{dt}\right|_{t=0}\mathcal D(B(t))=0.
	\]
	By the derivative formula for the logarithm of the determinant,
	\[
	\left.\frac{d}{dt}\right|_{t=0}\log\det M(t)
	=
	\operatorname{tr}\bigl(M^{-1}M'(0)\bigr),
	\]
	where
	\[
	M'(0)=\int_S\iota_EH_EP_E\,d\nu(E).
	\]
	Similarly,
	\[
	\left.\frac{d}{dt}\right|_{t=0}\log\det_E B_E(t)
	=
	\operatorname{tr}_E(B_E^{-1}H_E).
	\]
	The uniform positive-definiteness and boundedness of the fields justify
	differentiation under the integral. Consequently,
	\begin{equation}\label{eq:BB-first-variation}
		\int_S
		\operatorname{tr}_E\!\left[
		\bigl(P_EM^{-1}|_E-B_E^{-1}\bigr)H_E
		\right]d\nu(E)=0.
	\end{equation}
	Here $P_EM^{-1}|_E:E\to E$ denotes the compression of $M^{-1}$
	to $E$.
	
	Choose
	\[
	H_E=P_EM^{-1}|_E-B_E^{-1}.
	\]
	This is a continuous family of symmetric operators. By
	\eqref{eq:BB-first-variation},
	\[
	\int_S\operatorname{tr}_E(H_E^2)\,d\nu(E)=0.
	\]
	The integrand is nonnegative and vanishes if and only if $H_E=0$.
	It follows that $H_E=0$ for $\nu$-almost every $E$. Since $H_E$ is
	continuous and $S=\operatorname{supp}\nu$, we obtain
	\begin{equation}\label{eq:BB-compressed-inverse}
		B_E^{-1}=P_EM^{-1}|_E
		\qquad(E\in S).
	\end{equation}
	
	We next show that \eqref{eq:BB-compressed-inverse} implies that $M$
	commutes with $P_E$. Suppose first that $E\ne\R^n$. With respect to
	the orthogonal decomposition $\R^n=E\oplus E^\perp$, write
	\[
	M=
	\begin{pmatrix}
		C_E & R_E\\
		R_E^* & D'_E
	\end{pmatrix}.
	\]
	Here $C_E=P_EM|_E$, and $D'_E$ is positive definite on $E^\perp$.
	The Schur complement formula gives
	\[
	P_EM^{-1}|_E
	=
	\bigl(C_E-R_E(D'_E)^{-1}R_E^*\bigr)^{-1}.
	\]
	It follows from \eqref{eq:BB-compressed-inverse} that
	\begin{equation}\label{eq:BB-schur}
		B_E=C_E-R_E(D'_E)^{-1}R_E^*\le C_E.
	\end{equation}
	
	By the definition of $M$,
	\[
	\int_S\operatorname{tr}_E B_E\,d\nu(E)=\operatorname{tr}M.
	\]
	On the other hand, the isotropy condition gives
	\[
	\begin{aligned}
		\int_S\operatorname{tr}_E C_E\,d\nu(E)
		&=
		\int_S\operatorname{tr}(P_EM)\,d\nu(E)\\
		&=
		\operatorname{tr}\left(
		\left(\int_SP_E\,d\nu(E)\right)M
		\right)
		=\operatorname{tr}M.
	\end{aligned}
	\]
	Thus
	\[
	\int_S\operatorname{tr}_E(C_E-B_E)\,d\nu(E)=0.
	\]
	The operator $C_E-B_E$ is positive semidefinite by
	\eqref{eq:BB-schur}. Hence $C_E=B_E$ for $\nu$-almost every $E$,
	and then for every $E\in S$ by continuity. It follows that
	\[
	R_E(D'_E)^{-1}R_E^*=0.
	\]
	Since $(D'_E)^{-1}$ is positive definite, we obtain $R_E=0$. Therefore
	$M$ preserves both $E$ and $E^\perp$, and
	\[
	MP_E=P_EM,\qquad B_E=C_E=M|_E.
	\]
	If $E=\R^n$, the same conclusion follows directly from
	\eqref{eq:BB-compressed-inverse}. This proves the necessity of
	\eqref{eq:BB-equality}.
	
	Conversely, suppose that \eqref{eq:BB-equality} holds. Since $M$
	commutes with $P_E$ and $B_E=M|_E$, we have
	\[
	\log\det_E B_E=\operatorname{tr}(P_E\log M).
	\]
	Using the isotropy condition, we obtain
	\[
	\begin{aligned}
		\int_S\log\det_E B_E\,d\nu(E)
		&=
		\int_S\operatorname{tr}(P_E\log M)\,d\nu(E)\\
		&=
		\operatorname{tr}\left(
		\left(\int_SP_E\,d\nu(E)\right)\log M
		\right)\\
		&=\operatorname{tr}\log M
		=\log\det M.
	\end{aligned}
	\]
	Thus $\mathcal D(B)=0$, and equality holds in the Grassmannian
	Ball--Barthe inequality.
	
	We now extend the inequality to measurable positive definite fields.
	Assume that $E\mapsto B_E$ is measurable and that
	\[
	\int_S\operatorname{tr}_E B_E\,d\nu(E)<\infty,
	\qquad
	\int_S|\log\det_E B_E|\,d\nu(E)<\infty.
	\]
	
	We first consider the case in which there are constants $0<a<b<\infty$
	such that
	\[
	aI_E\le B_E\le bI_E
	\qquad\text{for $\nu$-almost every }E\in S.
	\]
	After changing the field on a null set, we may assume that these
	inequalities hold for every $E\in S$. Define the embedded logarithm
	\[
	L(E)=\iota_E(\log B_E)P_E.
	\]
	Then $L$ is a bounded measurable map from $S$ to the space of
	symmetric operators on $\R^n$.
	By Lusin's theorem, for each positive integer $m$ there is a compact
	set $K_m\subset S$ such that
	\[
	\nu(S\setminus K_m)<2^{-m}
	\]
	and $L$ is continuous on $K_m$. Extending the matrix entries from
	$K_m$ to $S$ and taking the symmetric part, we obtain continuous
	symmetric fields $L_m$ on $S$ such that
	\[
	L_m(E)=L(E)\quad(E\in K_m),
	\qquad
	\sup_{m,E}\|L_m(E)\|<\infty.
	\]
	For each $E\in S$, set
	\[
	B_{E,m}
	=
	\exp_E\bigl(P_EL_m(E)|_E\bigr).
	\]
	The field $B_{E,m}$ is continuous and positive definite. Moreover,
	there are constants $c,C>0$, independent of $m$ and $E$, such that
	\[
	cI_E\le B_{E,m}\le CI_E.
	\]
	Since $\sum_m\nu(S\setminus K_m)<\infty$, almost every $E$ belongs to
	$K_m$ for all sufficiently large $m$. For such $E$,
	\[
	B_{E,m}
	=
	\exp_E\bigl(P_EL(E)|_E\bigr)
	=B_E
	\]
	for all sufficiently large $m$.
	
	The inequality for continuous fields gives
	\[
	\det\left(
	\int_S\iota_EB_{E,m}P_E\,d\nu(E)
	\right)
	\ge
	\exp\left(
	\int_S\log\det_EB_{E,m}\,d\nu(E)
	\right).
	\]
	The uniform bounds above permit us to apply the dominated convergence
	theorem to both integrals. Letting $m\to\infty$, we obtain
	\[
	\det\left(
	\int_S\iota_EB_EP_E\,d\nu(E)
	\right)
	\ge
	\exp\left(
	\int_S\log\det_EB_E\,d\nu(E)
	\right)
	\]
	under the additional bounds $aI_E\le B_E\le bI_E$.
	
	For a general measurable positive definite field, let
	$\lambda_1(E),\ldots,\lambda_k(E)$ be the eigenvalues of $B_E$.
	For $r\ge1$, define $B_{E,r}$ by replacing each eigenvalue
	$\lambda_i(E)$ with
	\[
	\lambda_{i,r}(E)
	=
	\min\{r,\max\{r^{-1},\lambda_i(E)\}\}.
	\]
	Then
	\[
	r^{-1}I_E\le B_{E,r}\le rI_E,
	\]
	so the preceding case applies to $B_{E,r}$. Also,
	$B_{E,r}\to B_E$ for almost every $E$, and
	\[
	\operatorname{tr}_E B_{E,r}
	\le \operatorname{tr}_E B_E+k.
	\]
	Since $\nu(S)<\infty$, the dominated convergence theorem gives
	\[
	\int_S\iota_EB_{E,r}P_E\,d\nu(E)
	\longrightarrow
	\int_S\iota_EB_EP_E\,d\nu(E).
	\]
	
	Set
	\[
	P(E)=\sum_{i=1}^k\log^+\lambda_i(E),
	\qquad
	N(E)=\sum_{i=1}^k\log^-\lambda_i(E).
	\]
	Since $\log^+\lambda\le\lambda$ for $\lambda>0$,
	\[
	P(E)\le\operatorname{tr}_E B_E.
	\]
	Moreover,
	\[
	\log\det_EB_E=P(E)-N(E),
	\]
	and hence
	\[
	N(E)\le P(E)+|\log\det_EB_E|.
	\]
	Our assumptions therefore imply that
	\[
	\sum_{i=1}^k|\log\lambda_i(E)|=P(E)+N(E)
	\]
	is integrable on $S$.
	
	For each $i$, we have $\lambda_{i,r}(E)\to\lambda_i(E)$ and
	\[
	|\log\lambda_{i,r}(E)|
	\le |\log\lambda_i(E)|.
	\]
	Thus
	\[
	\log\det_EB_{E,r}\longrightarrow\log\det_EB_E
	\]
	almost everywhere, and
	\[
	|\log\det_EB_{E,r}|
	\le
	\sum_{i=1}^k|\log\lambda_i(E)|.
	\]
	Another application of the dominated convergence theorem yields
	\[
	\int_S\log\det_EB_{E,r}\,d\nu(E)
	\longrightarrow
	\int_S\log\det_EB_E\,d\nu(E).
	\]
	Applying the inequality to $B_{E,r}$ and letting $r\to\infty$
	completes the proof for measurable positive definite fields.
	\end{proof}

\begin{lemma}\label{lem:energy}
	For a square-integrable field $z_E\in E$ and $z=\int z_E\,d\nu$,
	\begin{equation}\label{eq:energy}
		\int|z_E|^2\,d\nu-|z|^2=\int|z_E-P_Ez|^2\,d\nu.
	\end{equation}
\end{lemma}
\begin{proof}
Since $z_E\in E$ and $z=\int_Sz_E\,d\nu(E)$,
\[
 \int_S\langle z_E,P_Ez\rangle\,d\nu(E)
 =\left\langle\int_Sz_E\,d\nu(E),z\right\rangle=|z|^2.
\]
By \eqref{eq:isotropy},
\[
 \int_S|P_Ez|^2\,d\nu(E)
 =\left\langle\left(\int_SP_E\,d\nu(E)\right)z,z\right\rangle=|z|^2.
\]
Expanding the square and using the two identities gives
\[
 \begin{aligned}
 \int_S|z_E-P_Ez|^2\,d\nu(E)
 &=\int_S|z_E|^2\,d\nu(E)
   -2\int_S\langle z_E,P_Ez\rangle\,d\nu(E)
   +\int_S|P_Ez|^2\,d\nu(E)\\
 &=\int_S|z_E|^2\,d\nu(E)-|z|^2.
 \end{aligned}
\]
This proves \eqref{eq:energy}.
\end{proof}


\begin{lemma}[Hessian rigidity]\label{lem:rigidity}
	Let
	\[
	\R^n=D\oplus F_1\oplus\cdots\oplus F_s
	\]
	be the orthogonal decomposition defined in
	\eqref{eq:decomposition}. For $x\in\R^n$, write
	\[
	x_D=P_Dx,\qquad x_j=P_{F_j}x,\quad j=1,\ldots,s.
	\]
	Suppose that $\Phi\in C^3(\R^n)$ and
	\[
	\nabla^2\Phi(x)P_E=P_E\nabla^2\Phi(x)
	\qquad(x\in\R^n,\ E\in S).
	\]
	Then there exist a constant $c\in\R$, a vector $d\in D$, a symmetric
	operator $C:D\to D$, and functions $\psi_j\in C^3(F_j)$ such that
	\begin{equation}\label{eq:rigidity}
		\Phi(x)
		=
		c+\frac12\langle Cx_D,x_D\rangle+\langle d,x_D\rangle
		+\sum_{j=1}^s\psi_j(x_j).
	\end{equation}
	The operator $C$ is independent of $x$ and satisfies
	\[
	CP_{E_D}=P_{E_D}C
	\qquad(E\in S),
	\]
	where $E_D=E\cap D$ and $P_{E_D}$ is the orthogonal projection in
	$D$ onto $E_D$. The conclusion does not require $\Phi$ to be convex.
\end{lemma}

\begin{proof}
Since $\nabla^2\Phi(x)$ commutes with every $P_E$, it preserves each complete common eigenspace $F_j$. As $\nabla^2\Phi(x)$ is symmetric, it also preserves $D$. Thus, for $i\ne j$,
 \begin{equation}\label{e-1}
 	P_D\nabla^2\Phi(x)P_{F_j}=0,
 \qquad P_{F_i}\nabla^2\Phi(x)P_{F_j}=0.
 \end{equation}

We next show that the restriction of $\nabla^2\Phi(x)$ to $D$ is independent of $x$. Fix $x\in\mathbb R^n$. For each $z\in\mathbb R^n$, define
\[
C_z=\partial_z(\nabla^2\Phi)(x)
=\left.\frac{d}{dt}\right|_{t=0}\nabla^2\Phi(x+tz),
\]
where $\partial_z$ denotes the directional derivative in the direction $z$.
Each $C_z$ is symmetric. Symmetry of the third derivatives and differentiation of the commutation relation give
\[
 C_zv=C_vz,\qquad C_zP_E=P_EC_z
 \quad(z,v\in\R^n,\ E\in S).
\]
Consequently,
\[
 C_{P_Ez}v=C_vP_Ez=P_EC_vz=P_EC_zv,
 \qquad C_{P_Ez}=P_EC_z.
\]
Let $E,L\in S$. Since $C_{P_Ez}$ commutes with $P_L$, we have
\[
C_{P_Ez}P_L=P_LC_{P_Ez}.
\]
Using $C_{P_Ez}=P_EC_z$ and $C_zP_L=P_LC_z$, we obtain
\[
P_EP_LC_z
=P_EC_zP_L
=C_{P_Ez}P_L
=P_LC_{P_Ez}
=P_LP_EC_z.
\]
Thus, for every $v\in\mathbb R^n$,
\[
P_EP_L(C_zv)=P_LP_E(C_zv).
\]
Let
\[
 W=\operatorname{span}\{C_zv:z,v\in\R^n\}.
\]
The identity $P_EC_z=C_{P_Ez}$ gives
\[
 P_E(C_zv)=C_{P_Ez}v\in W.
\]
Thus $P_EW\subset W$ for every $E$. The identity above also gives
\[
 P_EP_Lw=P_LP_Ew\qquad(w\in W,\ E,L\in S).
\]
The restrictions of the projections to $W$ are therefore commuting
symmetric operators. They have a common orthonormal eigenbasis.
Every eigenvalue of a projection is zero or one, so each vector in
this basis belongs to one of the complete spaces $F_j$. Hence
\[
 W\subset F_1\oplus\cdots\oplus F_s=D^\perp.
\]
Each symmetric $C_z$ has range in $D^\perp$. Thus, for $v\in D$
and arbitrary $u\in\R^n$,
\[
 \langle C_zv,u\rangle=\langle v,C_zu\rangle=0.
\]
It follows that $C_zv=0$ for every $v\in D$. In particular,
\[
 \partial_z\bigl(P_D\nabla^2\Phi(x)|_D\bigr)=0
 \qquad(x,z\in\R^n).
\]
There is therefore a constant symmetric operator $C:D\to D$ such that $\nabla^2\Phi(x)|_D=C$ for all $x$. Since $\nabla^2\Phi(x)$ commutes with the projections,
\[
 CP_{E_D}=P_{E_D}C\qquad(E\in S).
\]

Write
\[
x=x_D+x_1+\cdots+x_s,
\qquad
x_D\in D,\quad x_j\in F_j.
\]
For $1\le j\le s$, put
\[
a_j=x_D+x_1+\cdots+x_{j-1}.
\]
Identities \eqref{e-1} imply that
\[
P_{F_j}\nabla^2\Phi(y)a_j=0
\qquad (y\in\R^n).
\]
Hence, for every $t\in[0,1]$,
\[
\frac{d}{dr}
\left[P_{F_j}\nabla\Phi(ra_j+tx_j)\right]
=
P_{F_j}\nabla^2\Phi(ra_j+tx_j)a_j
=0.
\]
Taking $r=0$ and $r=1$, we obtain
\[
P_{F_j}\nabla\Phi(a_j+tx_j)
=
P_{F_j}\nabla\Phi(tx_j).
\]
Since $x_j\in F_j$, integration along the segment
$t\mapsto a_j+tx_j$ gives
\[
\begin{aligned}
	\Phi(a_j+x_j)-\Phi(a_j)
	&=\int_0^1
	\langle\nabla\Phi(a_j+tx_j),x_j\rangle\,dt\\
	&=\int_0^1
	\langle\nabla\Phi(tx_j),x_j\rangle\,dt\\
	&=\Phi(x_j)-\Phi(0).
\end{aligned}
\]
Summing these identities for $j=1,\ldots,s$, we find that
\[
\Phi(x)
=
\Phi(x_D)
+\sum_{j=1}^s\bigl(\Phi(x_j)-\Phi(0)\bigr).
\]

We now consider the restriction of $\Phi$ to $D$. Set
\[
c=\Phi(0),
\qquad
d=P_D\nabla\Phi(0)\in D.
\]
Since its Hessian is the constant operator $C$, Taylor's formula on
$D$ yields
\[
\Phi(x_D)
=
c+\langle d,x_D\rangle
+\frac12\langle Cx_D,x_D\rangle.
\]
For each $j$, define
\[
\psi_j(z)=\Phi(z)-\Phi(0),
\qquad z\in F_j.
\]
Then $\psi_j\in C^3(F_j)$, and
\[
\Phi(x)
=
c+\frac12\langle Cx_D,x_D\rangle
+\langle d,x_D\rangle
+\sum_{j=1}^s\psi_j(x_j).
\]
This proves \eqref{eq:rigidity}.
\end{proof}

We shall also use the following consequence of continuity. Suppose that a continuous probability density on an orthogonal product agrees almost everywhere with a product of measurable probability densities. Then each factor has a continuous version. By Fubini, fix the other coordinates where their factors are finite and positive and the factorization holds for almost every remaining coordinate. The resulting continuous slice supplies this version. The factorization then holds everywhere. An almost-everywhere identity between jointly continuous bundle expressions extends to all fibers, because every nonempty relatively open bundle neighborhood has positive parameter--fiber measure.

\begin{lemma}[Recovery of factors after Gaussian convolution]\label{lem:deconvolution}
	Let $V=V_0\oplus\cdots\oplus V_s$ be an orthogonal sum, and let $p$ be a probability density on $V$. For $x\in V$, write $x_j=P_{V_j}x$. Let $L$ be a Euclidean space of dimension $d$. For $t>0$, define
	the Gaussian density on $L$ by
	\[
	\gamma_{L,t}(z)
	=(4\pi t)^{-d/2}
	\exp\left(-\frac{|z|^2}{4t}\right),
	\qquad z\in L.
	\]  For $a\in L$ and a positive definite symmetric operator $\Sigma$ on $L$, write
\[
 \varphi_{L,a,\Sigma}(z)
 =\frac{\exp\!\left(-\frac12\langle\Sigma^{-1}(z-a),z-a\rangle\right)}
 {(2\pi)^{d/2}(\det_L\Sigma)^{1/2}},\qquad z\in L.
\]
Let $t>0$, and suppose
	\[
	p*\gamma_{V,t}(x)
	=\varphi_{V_0,a,\Sigma_*}(x_0)\prod_{j=1}^s g_j(x_j)
	\]
	where $a\in V_0$, $\Sigma_*>0$ on $V_0$, and each $g_j$ is a probability density on $V_j$.
	Then $\Sigma_0=\Sigma_*-2tI_{V_0}>0$ and
	\[
	p(x)=\varphi_{V_0,a,\Sigma_0}(x_0)\prod_{j=1}^s h_j(x_j)
	\quad\hbox{almost everywhere},
	\qquad g_j=h_j*\gamma_{V_j,t}.
	\]
	Here $h_j$ is the actual marginal density of $p$ on $V_j$. 
\end{lemma}
\begin{proof}

Write $\widehat p(\xi)=\int_V e^{i\langle\xi,x\rangle}p(x)\,dx$. Let $h_j$ be the marginal density of $p$ on $V_j$.	The Gaussian density $\gamma_{V,t}$ satisfies
	\[
	\widehat{\gamma_{V,t}}(\xi)=e^{-t|\xi|^2}.
	\]
	Taking Fourier transforms of the assumed convolution identity gives
	\[
	\widehat p(\xi)e^{-t|\xi|^2}
	=
	e^{i\langle a,\xi_0\rangle
		-\frac12\langle\Sigma_*\xi_0,\xi_0\rangle}
	\prod_{j=1}^s\widehat g_j(\xi_j).
	\]
	
	First, take $\xi=\xi_0\in V_0$. Since each $g_j$ is a probability
	density, $\widehat g_j(0)=1$. Hence
	\[
	\widehat p_0(\xi_0)e^{-t|\xi_0|^2}
	=
	e^{i\langle a,\xi_0\rangle
		-\frac12\langle\Sigma_*\xi_0,\xi_0\rangle}.
	\]
	With
	\[
	\Sigma_0=\Sigma_*-2tI_{V_0},
	\]
	this becomes
	\[
	\widehat p_0(\xi_0)
	=
	\exp\left(
	i\langle a,\xi_0\rangle
	-\frac12\langle\Sigma_0\xi_0,\xi_0\rangle
	\right).
	\]
	
	We show that $\Sigma_0$ is positive definite. Suppose first that
	$\Sigma_0$ has a negative eigenvalue $\lambda$, with corresponding
	unit eigenvector $v\in V_0$. Then, for $r\in\R$,
	\[
	\bigl|\widehat p_0(rv)\bigr|
	=
	\exp\left(-\frac12\lambda r^2\right).
	\]
	This is greater than $1$ when $r\ne0$, contradicting
	$|\widehat p_0|\le1$. Therefore
	\[
	\Sigma_0\ge0.
	\]
	
	Suppose next that $\Sigma_0$ has a nonzero kernel vector $v$.
	Let $X_0$ be a random vector in $V_0$ with density $p_0$. For every
	$r\in\R$,
	\[
	\begin{aligned}
		\mathbb E e^{ir\langle v,X_0-a\rangle}
		&=
		e^{-ir\langle v,a\rangle}\widehat p_0(rv)\\
		&=
		\exp\left(-\frac{r^2}{2}
		\langle\Sigma_0v,v\rangle\right)
		=1.
	\end{aligned}
	\]
	The right-hand side is the characteristic function of the constant
	random variable $0$. By uniqueness of characteristic functions,
	\[
	\langle v,X_0-a\rangle=0
	\qquad\text{almost surely}.
	\]
	Thus the law of $X_0$ is supported on the hyperplane
	\[
	\{u\in V_0:\langle v,u-a\rangle=0\}.
	\]
	This is impossible because $p_0$ is a probability density with
	respect to Lebesgue measure on $V_0$. Consequently,
	\[
	\Sigma_0>0.
	\]
	
	Now fix $j\ge1$ and take $\xi=\xi_j\in V_j$. The Fourier identity
	gives
	\[
	\widehat g_j(\xi_j)
	=
	\widehat p_j(\xi_j)e^{-t|\xi_j|^2}
	=
	\widehat h_j(\xi_j)e^{-t|\xi_j|^2}.
	\]
	By uniqueness of characteristic functions,
	\[
	g_j=h_j*\gamma_{V_j,t}
	\qquad\text{almost everywhere on }V_j.
	\]
	
	Finally, let $\xi=\xi_0+\cdots+\xi_s\in V$ be arbitrary. Substituting
	the identities for $\widehat g_j$ into the Fourier identity, we get
	\[
	\begin{aligned}
		\widehat p(\xi)e^{-t|\xi|^2}
		&=
		e^{i\langle a,\xi_0\rangle
			-\frac12\langle\Sigma_*\xi_0,\xi_0\rangle}
		\prod_{j=1}^s
		\left(\widehat h_j(\xi_j)e^{-t|\xi_j|^2}\right)\\
		&=
		e^{i\langle a,\xi_0\rangle
			-\frac12\langle\Sigma_0\xi_0,\xi_0\rangle}
		\prod_{j=1}^s\widehat h_j(\xi_j)
		\cdot e^{-t|\xi|^2}.
	\end{aligned}
	\]
		Since $e^{-t|\xi|^2}>0$, it follows that
	\[
	\widehat p(\xi)
	=
	e^{i\langle a,\xi_0\rangle
		-\frac12\langle\Sigma_0\xi_0,\xi_0\rangle}
	\prod_{j=1}^s\widehat h_j(\xi_j).
	\]
	The right-hand side is the characteristic function of the probability
	density
	\[
	x\longmapsto
	\varphi_{V_0,a,\Sigma_0}(x_0)
	\prod_{j=1}^s h_j(x_j).
	\]
	Uniqueness of characteristic functions therefore yields
	\[
	p(x)
	=
	\varphi_{V_0,a,\Sigma_0}(x_0)
	\prod_{j=1}^s h_j(x_j)
	\qquad\text{for almost every }x\in V.
	\]
	By construction, each $h_j$ is the marginal density of $p$ on $V_j$.
\end{proof}

\section{Parameter transport and measurable envelopes}\label{sec:transport}

In this section we prove the transport results needed below. For the forward inequality, we use maps between bounded balls. For the reverse inequality, we use convex transports with bounded targets and an area inequality for their average. We also study smooth transport maps for Gaussian-regularized families. These maps are used in the equality proofs.

\subsection{Measurability of the envelopes}
\begin{proposition}\label{prop:measurability}
	For an admissible family, both $\RB(f)$ and $\RBone(f)$ are universally measurable. The bounded-section assertion remains valid for a jointly Borel family that is continuous in the fiber variable and uniformly bounded on each fiber ball.
\end{proposition}
\begin{proof}
	We use the completion of $\nu$ and identify sections that agree
	$\nu$-almost everywhere. Every bounded measurable section has a Borel representative.
	If $b(E)$ is such a representative and $|v(E)|\le N$ almost
	everywhere, define
	\[
	 \widetilde v(E)=
	 \begin{cases}
	 P_Eb(E),&|P_Eb(E)|\le N,\\
	 0,&|P_Eb(E)|>N.
	 \end{cases}
	\]
	Then $\widetilde v(E)\in E$ and $|\widetilde v(E)|\le N$ for every
	$E$, while $\widetilde v=v$ almost everywhere. Thus the barycenter
	and logarithmic integral do not change.
	Let $X$ be the closed subspace of $L^2(\nu;\R^n)$ consisting of sections $v(E)\in E$. This is a separable Hilbert space. For each  integer $N$, set
	\[
	K_N=\left\{
	v\in X:\ |v(E)|\le N
	\text{ for }\nu\text{-almost every }E\in S
	\right\}.
	\]
	The set $K_N$ is closed in $X$ with respect to the $L^2$ norm.
	Since $X$ is a separable Hilbert space, $K_N$, equipped with the
	$L^2$ metric, is a Polish space. Since $f$ is continuous and
	\[
	\{(E,y):E\in S,\ y\in E,\ |y|\le N\}
	\]
	is compact, there is a constant $M_N\ge1$ such that
	\[
	f_E(y)\le M_N.
	\]
	Consequently, for every $v\in K_N$,
	\[
	J(v):=\int_S\log f_E(v(E))\,d\nu(E)
	\in[-\infty,\nu(S)\log M_N]
	\]
	 is upper semicontinuous. Indeed, any convergent sequence has a subsequence converging almost everywhere. Fiber continuity and the uniform upper bound permit reverse Fatou. Applying this argument to a subsequence realizing the limsup proves upper semicontinuity for the entire sequence. Define the barycenter map
\[
 B(v)=\int_Sv(E)\,d\nu(E)\in\R^n.
\]
For $u,v\in X$, the Cauchy--Schwarz inequality gives
\[
 |B(u)-B(v)|
 \le\int_S|u(E)-v(E)|\,d\nu(E)
 \le\nu(S)^{1/2}\|u-v\|_{L^2(\nu)}.
\]
Thus $B:X\to\R^n$ is continuous. The set
	 $
	 \{(x,v)\in\R^n\times K_N:\int v\,d\nu=x\}
	$
	 is closed. Moreover, the upper semicontinuity of $J$ implies that
	 \[
	 \{v\in K_N:J(v)\ge c\}
	 \]
	 is closed for every $c\in\R$. Hence, for every $a\in\R$,
	 \[
	 \{v\in K_N:J(v)>a\}
	 =
	 \bigcup_{m=1}^{\infty}
	 \{v\in K_N:J(v)\ge a+1/m\}
	 \]
	 is a Borel set. Therefore
	 \[
	 \{(x,v)\in\R^n\times K_N:B(v)=x,\ J(v)>a\}
	 \]
	 is also Borel. Write $\pi$ for projection onto the $x$ coordinate. For $a\in\R$,
	\[
	\{x:\RB(f)(x)>e^a\}
	=\bigcup_{N=1}^{\infty}\pi\left\{(x,v):v\in K_N,\ \int_Sv\,d\nu=x,\ J(v)>a\right\}.
	\]
	Each projected set is analytic and hence universally measurable. This proves the assertion for $\RB$.
	
	For $\RBone$, let $X_1$ be the closed subspace of
	$L^1(\nu;\R^n)$ consisting of sections. For $m\ge1$, define
	\[
	 J_m^\pm(v)=\int_S
	 \min\{\log^\pm f_E(v(E)),m\}\,d\nu(E),\qquad v\in X_1.
	\]
	If $v_i\to v$ in $L^1$, every subsequence has a further subsequence
	converging to $v$ almost everywhere. Fiber continuity and the bound
	$0\le\min\{\log^\pm f_E,m\}\le m$ show by dominated convergence
	that $J_m^\pm(v_i)\to J_m^\pm(v)$. Thus $J_m^\pm$ is continuous.
	By monotone convergence,
	\[
	 J_m^\pm(v)\uparrow
	 \int_S\log^\pm f_E(v(E))\,d\nu(E).
	\]
	The two extended-real limits are therefore Borel functions. In
	particular, the set on which the positive-part integral is finite is
	Borel. On this set $J$ is an extended-real Borel function. The map
	$v\mapsto B(v)$ is continuous on $X_1$, so the same projection
	argument proves the assertion for $\RBone$.

	For the last assertion, joint Borel measurability gives
	measurability of $E\mapsto f_E(v(E))$. Fiber continuity gives the
	pointwise convergence used above, and the uniform bound on each
	ball gives the common upper bound for reverse Fatou. Thus the same
	argument for $K_N$ proves the bounded-section assertion. 
\end{proof}

\subsection{Transport between fixed balls}
For a Euclidean space $L$ and $R>0$, write
\[
 B_L(R)=\{y\in L:|y|<R\},\qquad
 B^k(R)=\{y\in\R^k:|y|<R\}.
\]
\begin{lemma}\label{lem:ball-transport}
	Fix $R,\rho>0$. For each $E\in S$, let $a_E$ be a positive
	probability density on $\overline{B_E(R)}$. Assume that $a_E$
	depends continuously on $E$ in the $C^2$ norm in local orthonormal
	frames. Set
	\[
	Z_\rho
	=
	\int_{B^k(\rho)}e^{-\pi|z|^2}\,dz
	\]
	and define
	\[
	g_{E,\rho}(z)
	=
	Z_\rho^{-1}e^{-\pi|z|^2},
	\qquad z\in B_E(\rho).
	\]
	Let
	\[
	T_{E,\rho}:B_E(R)\longrightarrow B_E(\rho)
	\]
	be the convex-gradient transport from $a_E$ to $g_{E,\rho}$.
	Then $T_{E,\rho}$ is a diffeomorphism with positive definite
	derivative. The maps $T_{E,\rho}$ and their first derivatives in
	the fiber variables are jointly continuous in local orthonormal
	frames. Moreover,
	\begin{equation}\label{eq:ball-Jacobian}
		a_E(y)
		=
		Z_\rho^{-1}e^{-\pi|T_{E,\rho}(y)|^2}
		\det_E DT_{E,\rho}(y),
		\qquad y\in B_E(R).
	\end{equation}
\end{lemma}

\begin{proof}
Fix $E_0\in S$. In a local orthonormal frame near $E_0$, the
source and target domains are the fixed Euclidean balls
$B^k(R)$ and $B^k(\rho)$. We work in these coordinates and write
\[
g_\rho(z)=Z_\rho^{-1}e^{-\pi|z|^2},
\qquad z\in B^k(\rho).
\]

Let $\phi_{E,\rho}$ be the convex transport potential from
$a_E$ to $g_\rho$, normalized by
\[
\phi_{E,\rho}(0)=0,
\qquad
T_{E,\rho}=\nabla\phi_{E,\rho}.
\]
The transport equation is
\[
\det \nabla^2\phi_{E,\rho}(x)
=
\frac{a_E(x)}
{g_\rho(\nabla\phi_{E,\rho}(x))}
\]
with boundary condition
\[
\nabla\phi_{E,\rho}(B^k(R))=B^k(\rho).
\]

The two balls are uniformly convex. Moreover, the  densities $a_E$
have uniform positive lower bounds and uniform $C^2$ bounds, while
$g_\rho$ is fixed, smooth and positive on
$\overline{B^k(\rho)}$. By the global estimates for the second
boundary value problem \cite{Caffarelli,TW}, there are constants
$\alpha\in(0,1)$ and $C>0$, independent of $E$, such that
\[
\|\phi_{E,\rho}\|_{C^{2,\alpha}(\overline{B^k(R)})}\le C.
\]

We next show that the potentials depend continuously on $E$.
Suppose that $E_i\to E$ in this local frame. The uniform estimate
and compactness give a subsequence, still denoted by $E_i$, such
that
\[
\phi_{E_i,\rho}\longrightarrow\phi
\qquad\text{in }C^2(\overline{B^k(R)})
\]
for some convex function $\phi$.

For every continuous function $\eta$ on $\overline{B^k(\rho)}$,
the transport identity gives
\[
\int_{B^k(R)}
\eta\bigl(\nabla\phi_{E_i,\rho}(z)\bigr)a_{E_i}(z)\,dz
=
\int_{B^k(\rho)}\eta(y)g_\rho(y)\,dy.
\]
Since $a_{E_i}\to a_E$ and
$\nabla\phi_{E_i,\rho}\to\nabla\phi$ uniformly, we may pass to
the limit to obtain
\[
\int_{B^k(R)}
\eta(\nabla\phi(z))a_E(z)\,dz
=
\int_{B^k(\rho)}\eta(y)g_\rho(y)\,dy.
\]
Thus $\nabla\phi$ transports $a_E$ to $g_\rho$. By uniqueness
of the normalized convex transport potential,
\[
\phi=\phi_{E,\rho}.
\]
Every convergent subsequence has the same limit. Therefore
\[
\phi_{E_i,\rho}\longrightarrow\phi_{E,\rho}
\qquad\text{in }C^2(\overline{B^k(R)}).
\]
In particular, $T_{E,\rho}$ and $DT_{E,\rho}$ depend continuously
on $E$ in local orthonormal frames.

The Jacobian equation gives
\[
 \det DT_{E,\rho}(x)
 =\frac{a_E(x)}{g_\rho(T_{E,\rho}(x))}>0.
\]
Since $DT_{E,\rho}$ is symmetric and nonnegative by convexity,
it is positive definite. Orthogonal changes of frame preserve
$g_\rho$, so uniqueness makes the local constructions agree.

For $k=1$ the same conclusions follow directly. With an oriented
coordinate on $E$, define
\[
 F_E(x)=\int_{-R}^x a_E(s)\,ds,\qquad
 G(x)=\int_{-\rho}^x g_\rho(s)\,ds.
\]
Then $T_{E,\rho}=G^{-1}\circ F_E$ and
\[
 T_{E,\rho}'(x)=\frac{a_E(x)}{g_\rho(T_{E,\rho}(x))}>0.
\]
The assumed $C^2$ convergence of the source densities gives convergence
of these maps and their first derivatives. The result is unchanged if
the orientation is reversed.
\end{proof}

\subsection{Weak transport with bounded targets}
We shall use families $q=(q_E)_{E\in S}$ of probability densities
satisfying the following conditions. The map $(E,y)\mapsto q_E(y)$ is
Borel measurable on $\mathcal T$, each $q_E$ is continuous on $E$, and
\[
 \int_Eq_E(y)\,dy=1\qquad(E\in S).
\]
For every $R>0$ there is a constant $C_R$ such that
\[
 q_E(y)\le C_R\qquad(E\in S,\ y\in E,\ |y|\le R).
\]
Finally, the family is uniformly tight in the sense of
\eqref{eq:tightness}. Continuity with respect to $E$ and moment
assumptions are not required.

\begin{lemma}\label{lem:weak-transport}
	Suppose $(E,y)\mapsto a_E(y)$ is a nonnegative Borel family of
	probability densities on $E$. Assume that there are $R,K>0$ such that
	\[
	 a_E(y)=0\quad\text{for }|y|>R,\qquad
	 a_E(y)\le K\quad(E\in S,\ y\in E).
	\] There are jointly measurable convex potentials $\phi_E$, normalized by $\phi_E(0)=0$, whose gradients $\nabla\phi_E$ transport $g_E(z)=e^{-\pi|z|^2}$ to $a_E$. They are uniformly Lipschitz. For almost every $(E,x)$, the Alexandrov Hessians $\nabla^2\phi_E(P_Ex)$ are positive definite and
	\begin{equation}\label{eq:weak-Jacobian}
		a_E(\nabla\phi_E(P_Ex))\det_E\nabla^2\phi_E(P_Ex)=e^{-\pi|P_Ex|^2}.
	\end{equation}
	For $\Phi(x)=\int_S\phi_E(P_Ex)\,d\nu$, almost everywhere
	\begin{equation}\label{eq:weak-Hessian}
		\nabla\Phi(x)=\int_S\iota_E\nabla\phi_E(P_Ex)\,d\nu,\qquad
		\nabla^2_a\Phi(x)\ \succeq\ \int_S\iota_E\nabla^2\phi_E(P_Ex)P_E\,d\nu.
	\end{equation}
	Here $\nabla^2_a$ denotes the density of the absolutely continuous part of the distributional Hessian.
\end{lemma}
\begin{proof}
Let
	\[
	d\gamma(x)=e^{-\pi|x|^2}\,dx
	\qquad(x\in\R^k)
	\]
and 	let $\mu$ be a probability measure supported in
	$\overline{B^k(R)}$. Both $\gamma$ and $\mu$ have finite second
	moments. 
	By Brenier's theorem \cite{Brenier}, there is a convex
	function $\phi$, normalized by $\phi(0)=0$, such that  $$(\nabla\phi)_\#\gamma=\mu.$$  (The pushforward measure $(\nabla\phi)_\#\gamma$ is defined by
$
	((\nabla\phi)_\#\gamma)(A)=\gamma((\nabla\phi)^{-1}(A))
$
	for every Borel set $A\subset \overline{B^k(R)}$.) 
	Since $\mu$ is supported in $\overline{B^k(R)}$, the transport
	identity gives
	\[
	|\nabla\phi(x)|\le R
	\qquad\text{for }\gamma\text{-almost every }x.
	\]
 By convexity, the potential
	can be chosen as an $R$-Lipschitz function on $\R^k$,  satisfying
	\[
	|\phi(x)-\phi(x')|\le R|x-x'|,
	\qquad
	|\phi(x)|\le R|x|
	\quad(x,x'\in\R^k).
	\]
	
 Let $\mu_i$ and $\mu$ be probability measures
	supported in $\overline{B^k(R)}$, and suppose that
	\[
	\mu_i\rightharpoonup\mu.
	\]
	Let $\phi_i$ and $\phi$ be the normalized convex potentials
	satisfying
	\[
	(\nabla\phi_i)_\#\gamma=\mu_i,
	\qquad
	(\nabla\phi)_\#\gamma=\mu,
	\qquad
	\phi_i(0)=\phi(0)=0.
	\]
	We claim that
	\[
	\phi_i\longrightarrow\phi
	\qquad\text{locally uniformly on }\R^k.
	\]
	
	Indeed, every $\phi_i$ is $R$-Lipschitz and vanishes at the origin.
	Hence, on each bounded set, the functions $\phi_i$ are uniformly
	bounded and equicontinuous. By the Arzel\`a--Ascoli theorem and a
	diagonal argument, every subsequence has a further subsequence,
	still denoted by $\phi_i$, such that
	\[
	\phi_i\longrightarrow\psi
	\qquad\text{locally uniformly on }\R^k,
	\]
	where $\psi$ is convex and $\psi(0)=0$.	
	For each $i$, define the transport plan
	\[
	\pi_i=(\operatorname{id},\nabla\phi_i)_\#\gamma.
	\]
	Its first and second marginals are $\gamma$ and $\mu_i$,
	respectively. Since $\gamma$ is fixed and every $\mu_i$ is supported
	in $\overline{B^k(R)}$, the plans $\pi_i$ form a tight family.
	Passing to a further subsequence if necessary, we have
	\[
	\pi_i\rightharpoonup\pi.
	\]
	The two marginals of $\pi$ are $\gamma$ and $\mu$.
	For $\pi_i$-almost every $(x,y)$, we have
	$y=\nabla\phi_i(x)$. The subgradient inequality therefore gives
	\[
	\phi_i(z)\ge\phi_i(x)+\langle y,z-x\rangle
	\qquad(z\in\R^k).
	\]
	Using the local uniform convergence of $\phi_i$ and the weak
	convergence of $\pi_i$, we may pass to the limit in this inequality and obtain
	\[
	\psi(z)\ge\psi(x)+\langle y,z-x\rangle
	\quad\text{for every }z\in\R^k
	\]
	for $\pi$-almost every $(x,y)$. A convex function is differentiable
	Lebesgue almost everywhere. Since the first marginal of $\pi$ is
	$\gamma$ and $\gamma$ has a density, $\psi$ is differentiable at
	$x$ for $\pi$-almost every $(x,y)$. At such a pair, for any
	$h\in\R^k$ and $t>0$, choose $z=x\pm th$ in the inequality above.
	Then
	\[
	 \frac{\psi(x)-\psi(x-th)}{t}\le\langle y,h\rangle
	 \le\frac{\psi(x+th)-\psi(x)}{t}.
	\]
	Letting $t\downarrow0$ gives
	$\langle y,h\rangle=\langle\nabla\psi(x),h\rangle$ for every $h$.
	Thus
	\[
	 y=\nabla\psi(x)\qquad\text{for }\pi\text{-almost every }(x,y).
	\]
	The second marginal of $\pi$ is $\mu$. Consequently,
	\[
	(\nabla\psi)_\#\gamma=\mu.
	\]
	By uniqueness of the Brenier transport from $\gamma$ to $\mu$,
	\[
	\nabla\psi=\nabla\phi
	\qquad\gamma\text{-almost everywhere}.
	\]
	The density of $\gamma$ is positive everywhere, so the equality also
	holds Lebesgue almost everywhere. Since $\psi$ and $\phi$ are locally
	Lipschitz, they differ by a constant. Their values at the origin are
	both zero. Hence
	\[
	\psi=\phi.
	\]
	Thus every subsequence of $(\phi_i)$ has a further subsequence that
	converges locally uniformly to $\phi$. It follows that the whole
	sequence satisfies
	\[
	\phi_i\longrightarrow\phi
	\qquad\text{locally uniformly on }\R^k.
	\]
	This argument also applies when the target measures do not have
	densities.

We next choose the potentials and their derivatives jointly
measurably in $E$ and the fiber variable.
Fix a local orthonormal frame and let
\[
U_E:\R^k\longrightarrow E
\]
be the corresponding isometry. In these coordinates, the target
measure has density
\[
\widetilde a_E(z)=a_E(U_Ez).
\]
Since $(E,z)\mapsto\widetilde a_E(z)$ is Borel measurable and the
target measures are supported in a fixed ball, the map
\[
E\longmapsto\widetilde\mu_E
:=\widetilde a_E(z)\,dz
\]
is Borel measurable for the weak topology.
Let $\widetilde\phi_E$ be the normalized convex potential from
$\gamma$ to $\widetilde\mu_E$. By the continuous dependence on the
target measure proved above,
\[
E\longmapsto\widetilde\phi_E
\]
is Borel measurable.
Consequently,
\[
(E,z)\longmapsto\widetilde\phi_E(z)
\]
is jointly Borel measurable.

We first choose a Borel representative of the gradients. For
$j=1,\ldots,k$ and $m\ge1$, define
\[
G_{j,m}(E,z)
=
m\bigl[\widetilde\phi_E(z+e_j/m)-\widetilde\phi_E(z)\bigr].
\]
Each $G_{j,m}$ is Borel measurable. Let $G_j(E,z)$ be the limit
of $G_{j,m}(E,z)$ when this limit exists and is finite, and set
$G_j(E,z)=0$ otherwise. Then $G_j$ is Borel measurable. At every
differentiability point of $\widetilde\phi_E$,
\[
G_j(E,z)=\frac{\partial\widetilde\phi_E}{\partial z_j}(z).
\]
For the Alexandrov Hessians, put
\[
Q_{u,m}(E,z)
=
m^2\bigl[
\widetilde\phi_E(z+u/m)
+\widetilde\phi_E(z-u/m)
-2\widetilde\phi_E(z)
\bigr].
\]
Each $Q_{u,m}$ is Borel measurable. Define a symmetric matrix field
$H(E,z)$ by
\[
H_{ii}(E,z)=\lim_{m\to\infty}Q_{e_i,m}(E,z)
\]
and, for $i\ne j$,
\[
H_{ij}(E,z)
=
\frac12\lim_{m\to\infty}
\bigl(
Q_{e_i+e_j,m}(E,z)
-Q_{e_i,m}(E,z)
-Q_{e_j,m}(E,z)
\bigr),
\]
whenever the corresponding limits exist and are finite. Set an
entry equal to zero when its defining limit does not exist as a
finite number. Thus $H$ is jointly Borel measurable.

At every Alexandrov second differentiability point of
$\widetilde\phi_E$, the second-order expansion gives
\[
\lim_{m\to\infty}Q_{u,m}(E,z)
=
\langle\nabla^2\widetilde\phi_E(z)u,u\rangle.
\]
Hence, by polarization,
\[
H(E,z)=\nabla^2\widetilde\phi_E(z)
\]
at such points. For each fixed $E$, convexity gives
differentiability almost everywhere, and Alexandrov's theorem gives
second differentiability almost everywhere. Therefore $G$ and $H$
represent the gradient and the Alexandrov Hessian almost everywhere
in each fiber.

Cover $S$ by finitely many frame neighborhoods and choose a
Borel partition subordinate to this cover. On each member of the
partition, use the corresponding coordinate construction and then
return to the subspace $E$. This gives jointly Borel representatives
of
\[
\phi_E(y),\qquad \nabla\phi_E(y),\qquad \nabla^2\phi_E(y).
\]
On overlaps, uniqueness of the normalized transport potential makes
the potentials agree after the orthogonal change of coordinates.
The derivative representatives agree at points where the
corresponding derivatives exist.
Their values at other points are used only to define measurable
representatives. The transport and Jacobian identities are asserted
only at the points where those identities have been established.

We next prove the Jacobian identity for a fixed $E\in S$. Write
\[
 d\gamma_E(z)=g_E(z)\,dz,\qquad d\mu_E(y)=a_E(y)\,dy,
 \qquad T_E=\nabla\phi_E.
\]
The transpose of the optimal plan induced by $T_E$ is optimal from
$\mu_E$ to $\gamma_E$. Since $\mu_E$ is absolutely continuous,
Brenier's theorem gives a reverse map $S_E$ and
\[
 S_E(T_E(z))=z\qquad\text{for }\gamma_E\text{-almost every }z.
\]
Consequently $T_E$ is injective on a Borel set of full $\gamma_E$-measure.

Take a Borel selection of the subgradient of $\phi_E$; it agrees with
$T_E$ at differentiability points. Alexandrov's theorem gives its
derivative as the Alexandrov Hessian almost everywhere. The area
formula on the essential injectivity set and the density formula for
the pushforward measure \cite[Lemma~1.3 and Theorem~1.4]{AS} imply
\[
 \det_E\nabla^2\phi_E(z)>0,\qquad
 a_E(T_E(z))\det_E\nabla^2\phi_E(z)=g_E(z)
 \quad\text{for almost every }z\in E.
\]
Here Lebesgue and $\gamma_E$ have the same null sets because $g_E>0$.
Since $\phi_E$ is convex, its Alexandrov Hessian is positive
semidefinite; the positive determinant therefore implies
\[
 \nabla^2\phi_E(z)>0\qquad\text{for almost every }z\in E.
\]
No convexity of the support of $\mu_E$ is used.

We now make the exceptional sets measurable in both variables. Let
$T_E(y)$ and $H_E(y)$ be the jointly Borel representatives of the
gradient and Alexandrov Hessian constructed above. In a local
orthonormal frame, the second-order expansion
\[
 \phi_E(y+h)=\phi_E(y)+\langle T_E(y),h\rangle
 +\tfrac12\langle H_E(y)h,h\rangle+o(|h|^2)
\]
defines a Borel set of pairs $(E,y)$. Indeed, the remainder condition
can be expressed using countably many rational vectors $h$ in the
frame. Using a finite Borel partition subordinate to the frame cover,
we obtain a Borel set $\mathcal N\subset\mathcal T$ on whose complement
the expansion holds, $H_E(y)>0$, and
\[
 a_E(T_E(y))\det_E H_E(y)=e^{-\pi|y|^2}.
\]
The last two conditions are Borel because the chosen representatives
and $a_E$ are jointly Borel. By Alexandrov's theorem and the fixed-fiber
identity, each fiber $N_E=\{y\in E:(E,y)\in\mathcal N\}$ is null.

For a fixed $E$, decompose $x=y+w\in E\oplus E^\perp$. Since
$P_Ex=y$ and $|N_E|_E=0$, Fubini's theorem gives
\[
 \left|\{x\in B^n(m):P_Ex\in N_E\}\right|=0
 \qquad(m\ge1).
\]
The set $\{(E,x):(E,P_Ex)\in\mathcal N\}$ is Borel. A second
application of Fubini's theorem therefore yields
\[
 \nu\{E\in S:P_Ex\in N_E\}=0
 \qquad\text{for almost every }x\in\R^n.
\]
In particular, for almost every $x$, both
\[
 H_E(P_Ex)>0,\qquad
 a_E(T_E(P_Ex))\det_EH_E(P_Ex)=e^{-\pi|P_Ex|^2}
\]
hold for $\nu$-almost every $E$.

Set
\[
\Phi(x)=\int_S\phi_E(P_Ex)\,d\nu(E),\qquad x\in\R^n.
\]
Each $\phi_E\circ P_E$ is convex and $R$-Lipschitz. Thus we may differentiate
under the parameter integral in the distributional sense. Using the
jointly measurable representatives chosen above, we obtain
\[
\nabla\Phi(x)
=
\int_S\iota_E \nabla\phi_E(P_Ex)\,d\nu(E)
\qquad\text{for almost every }x\in\R^n.
\]

Let $\mathbf H_E$ denote the distributional Hessian of $\phi_E\circ P_E$.
It is a nonnegative symmetric matrix-valued measure. For every
nonnegative compactly supported smooth function $\eta$,
\[
\begin{aligned}
	\int_{\R^n}\eta\,d\operatorname{tr}\mathbf H_E
	=
	-\int_{\R^n}\nabla\eta\cdot\nabla (\phi_E\circ P_E)\,dx\le R\int_{\R^n}|\nabla\eta|\,dx.
\end{aligned}
\]
Hence the measures $\mathbf H_E$ have uniformly bounded mass on
each compact set, independently of $E$.
For a convex function, the Alexandrov Hessian is the density of the
absolutely continuous part of its distributional Hessian; see
\cite[Section~3.2]{DPF}. Consequently, for every $\xi\in\R^n$,
\[
\langle\mathbf H_E\xi,\xi\rangle
\ge
\langle \iota_E
\nabla^2\phi_E(P_Ex)P_E\xi,\xi\rangle\,dx.
\]
Integrating over $E$ and using distributional differentiation under
the integral, we obtain
\[
\langle\mathbf H_\Phi\xi,\xi\rangle
\ge
\left(\int_S\langle \iota_E
\nabla^2\phi_E(P_Ex)P_E\xi,\xi\rangle\,d\nu(E)\right)dx,
\]
where $\mathbf H_\Phi$ is the distributional Hessian of $\Phi$.
For a compact set $K$, choose $\eta\in C_c^\infty(\R^n)$ with
$\eta\ge1$ on $K$ and $\eta\ge0$. Taking the trace in the preceding
measure inequality and using the uniform estimate gives
\[
 \int_K\int_S\operatorname{tr}_E\nabla^2\phi_E(P_Ex)\,d\nu(E)\,dx
 \le R\nu(S)\int_{\R^n}|\nabla\eta|\,dx<\infty.
\]
Thus the positive matrix integral is finite for almost every $x$.
Taking the densities of the absolutely continuous parts gives
\[
\langle\nabla_a^2\Phi(x)\xi,\xi\rangle
\ge
\int_S
\left\langle
\iota_E\nabla^2\phi_E(P_Ex)P_E\xi,P_E\xi
\right\rangle\,d\nu(E)
\]
for almost every $x$. Taking first a countable dense set of
$\xi\in\R^n$ and then using continuity of the quadratic forms, we
obtain
\[
\nabla_a^2\Phi(x)
\succeq
\int_S
\iota_E\nabla^2\phi_E(P_Ex)P_E\,d\nu(E)
\qquad\text{for almost every }x\in\R^n.
\]

The inequality may be strict. Averaging the singular parts of the
measures $\mathbf H_E$ can produce an absolutely continuous part.
\end{proof}

\begin{lemma}[Area inequality for the averaged convex gradient]\label{lem:convex-area}
	Let $\Phi:\R^n\to\R$ be finite and convex, with $\nabla^2_a\Phi>0$ almost everywhere. For every nonnegative universally measurable $H$,
	\begin{equation}\label{eq:convex-area}
		\int_{\R^n}H(y)\,dy\ge
		\int_{\R^n}H(\nabla\Phi(x))\det \nabla^2_a\Phi(x)\,dx
	\end{equation}
	holds, with the integrand on the right defined arbitrarily on the null exceptional set.
\end{lemma}
\begin{proof}
By Alexandrov's theorem and the assumption
$\nabla_a^2\Phi>0$ almost everywhere, we may choose a Borel set
$G\subset\R^n$ of full Lebesgue measure such that, for every
$x\in G$, the function $\Phi$ is twice differentiable in the
Alexandrov sense and
$
\nabla_a^2\Phi(x)>0.
$
Moreover, at each such point, the gradient is differentiable
relative to the set where it is defined, with derivative
$\nabla_a^2\Phi(x)$; see \cite[Section~1.3]{AS}. We show that the gradient is injective on $G$. Suppose that
$x,y\in G$ and $\nabla\Phi(x)=\nabla\Phi(y)$. The two supporting-plane
inequalities give
\[
 \Phi(y)\ge\Phi(x)+\langle\nabla\Phi(x),y-x\rangle,
 \qquad
 \Phi(x)\ge\Phi(y)+\langle\nabla\Phi(y),x-y\rangle.
\]
They must both be equalities. By convexity, $\Phi$ is therefore
affine on $[x,y]$. For $0<s<1$ its second-order expansion at $x$ would then give
	\[
	0=\tfrac12s^2\langle \nabla^2_a\Phi(x)(y-x),y-x\rangle+o(s^2),
	\]
	a contradiction.
	
	Decompose $G$, after discarding a null set if necessary, into disjoint Borel sets $G_\ell$ on which the gradient is Lipschitz. 
	Such a decomposition follows from differentiability: stratify the bound for the derivative and the radius in its defining remainder estimate, and intersect with cubes of sufficiently small diameter. Lipschitz extensions have derivative $\nabla^2_a\Phi$ at almost every density point of each $G_\ell$. The area formula on each piece, and injectivity across different pieces, therefore give the measure identity
	\begin{equation}\label{eq:area-pushforward}
		(\nabla\Phi)_\#\big({\bf1}_{G}\det \nabla^2_a\Phi\,dx\big)
		={\bf1}_{\nabla\Phi(G)}\,dy .
	\end{equation}
	Here $G$ denotes the union of the pieces. Its image is Borel by the injective Borel mapping theorem. For every nonnegative Borel function $H$, the identity gives
	\[
	\int_G H(\nabla\Phi(x))\det \nabla^2_a\Phi(x)\,dx
	=\int_{\nabla\Phi(G)}H(y)\,dy
	\le\int_{\R^n}H(y)\,dy.
	\]
	
	The image measure in \eqref{eq:area-pushforward} is dominated by Lebesgue measure. Hence the identity extends to the completions of these measures and to nonnegative universally measurable functions. This proves \eqref{eq:convex-area}. 
\end{proof}

\begin{remark}
	The area inequality can be strict even when the absolutely continuous Hessian is everywhere positive definite on its domain of definition. In one dimension $\Phi(x)=x^2/2+|x|$ has $\nabla^2_a\Phi=1$ almost everywhere, but the gradient on its differentiability set misses $[-1,1]$. The right side of \eqref{eq:convex-area} can therefore be strictly smaller than the left side. Only the inequality is needed below.
\end{remark}

\subsection{Smooth transport for heat-regularized families}
For $t>0$, let
\[
\gamma_{d,t}(z)
=(4\pi t)^{-d/2}e^{-|z|^2/(4t)}
\]
and, for $E\in S$, define
\[
q_{E,t}(y)
=(q_E*\gamma_{k,t})(y)
=\int_Eq_E(z)\gamma_{k,t}(y-z)\,dz,
\qquad y\in E.
\]
Then $q_{E,t}$ is a positive smooth probability density on $E$.
Fix $t>0$. In a local
orthonormal frame, the normalized densities $q_E$ are continuous in
$L^1$ by Scheff\'e's lemma. Thus, if $E_i\to E$,
\[
\|q_{E_i,t}-q_{E,t}\|_1
\le
\|q_{E_i}-q_E\|_1\|\gamma_{k,t}\|_1
=
\|q_{E_i}-q_E\|_1
\longrightarrow0.
\]
Moreover, for every multi-index $\alpha$,
\[
\|\partial^\alpha q_{E_i,t}
-\partial^\alpha q_{E,t}\|_\infty
\le
\|q_{E_i}-q_E\|_1
\|\partial^\alpha\gamma_{k,t}\|_\infty
\longrightarrow0.
\]
Hence $q_{E,t}$ depends continuously on $E$ in $L^1$ and in every
local $C^r$ norm.
If
$y=z+w$ and $|y|>R$, then $|z|>R/2$ or $|w|>R/2$. Therefore
\[
\int_{\{|y|>R\}\cap E}q_{E,t}(y)\,dy
\le
\int_{\{|z|>R/2\}\cap E}q_E(z)\,dz
+
\int_{\{|w|>R/2\}\cap E}\gamma_{k,t}(w)\,dw.
\]
By \eqref{eq:tightness}, the first term tends to zero uniformly in
$E$ as $R\to\infty$. The second term also tends to zero and is
independent of $E$ after identifying $E$ with $\R^k$. Thus
\[
\lim_{R\to\infty}
\sup_{E\in S}
\int_{\{|y|>R\}\cap E}q_{E,t}(y)\,dy=0.
\]
By \eqref{eq:tightness} again, we may choose $R_0>0$ such that
\[
\int_{B_E(R_0)}q_E(z)\,dz\ge\frac12
\qquad(E\in S).
\]
For every $M>0$, if $y\in E$ and $|y|\le M$, then
\[
\begin{aligned}
	q_{E,t}(y)
	&\ge
	\int_{B_E(R_0)}q_E(z)\gamma_{k,t}(y-z)\,dz\\
	&\ge
	\frac12(4\pi t)^{-k/2}
	\exp\left(-\frac{(M+R_0)^2}{4t}\right)>0.
\end{aligned}
\]
Thus $q_{E,t}$ has a positive lower bound on
$\{y\in E:|y|\le M\}$, uniform in $E\in S$.

\begin{lemma}\label{lem:heat-transport}
	For each fixed $t>0$, the Brenier--McCann maps $U_{E,t}:E\to E$ from $e^{-\pi|z|^2}\,dz$ to $q_{E,t}(y)\,dy$ are smooth diffeomorphisms. In local frames their derivatives of every fixed order depend continuously on the parameter on compact sets. The same is true of their inverses.
\end{lemma}
\begin{proof}
	Fix $t>0$ and let $E_i\to E$ in $S$.
	 Existence and uniqueness follow from McCann \cite{McCann}. A positive continuous density $g$ bounded above and below on a
ball $B$ is locally doubling on ellipsoids. That is,  if an ellipsoid $A$ and
its half dilation $\tfrac12 A$ lie in $B$, then
\[
 \frac{\int_A g}{\int_{\frac12 A}g}
 \le 2^k\frac{\sup_Bg}{\inf_Bg}.
\]
The source and target densities are smooth and positive, so they
satisfy the hypotheses of Theorems 1.6--1.7 of \cite{FJ}. Those
theorems, with both domains equal to $\R^k$, give smooth
diffeomorphisms $U_i=\nabla\phi_i$ onto $\R^k$. The same holds for
the limiting density. Normalize $\phi_i(0)=0$.
	
	\emph{Uniform bounds on compact source sets.}
	Fix $x$ with $|x|\le R$. Let $
	e=\frac{U_i(x)}{|U_i(x)|}
	$ with $U_i(x)\neq0$. 	 Monotonicity gives, for $z\in B(x+e,1/4)$,
	\[
	U_i(z)\cdot(z-x)\ge U_i(x)\cdot(z-x)\ge\tfrac34|U_i(x)|,
	\qquad |z-x|\le\tfrac54 .
	\]
	Consequently $|U_i(z)|\ge3|U_i(x)|/5$. The Gaussian mass of this ball is bounded below uniformly in $i$ and $|x|\le R$:
	\[
	\kappa_R=|B^k(1/4)|\min_{|z|\le R+2}e^{-\pi|z|^2}>0.
	\]
	Uniform tightness of the target family provides $L_R$ with
	$\sup_i\int_{|y|>L_R}q_{E_i,t}(y)\,dy<\kappa_R$.
	If $|U_i(x)|>5L_R/3$, the pushforward identity gives
	\[
	\kappa_R\le\int_{B(x+e,1/4)}e^{-\pi|z|^2}\,dz
	\le\int_{|y|>L_R}q_{E_i,t}(y)\,dy<\kappa_R,
	\]
	a contradiction. Hence $|U_i(x)|\le5L_R/3$ for $|x|\le R$. Since $\phi_i(0)=0$, the potentials are locally uniformly bounded as well. 
	
	\emph{Identification of the potential limit.}
	Every subsequence has a further subsequence whose normalized
potentials converge locally uniformly. Its transport plans are tight
because their two marginals are tight. The argument used in the proof
of Lemma~\ref{lem:weak-transport} passes the supporting-plane
inequality to a weak limit of these plans. The limit has the Gaussian
source and target $q_{E,t}$ and is induced by the gradient of the
limiting potential. McCann uniqueness identifies it with the map for
$q_{E,t}$. Since the normalized potential is unique, the whole
sequence converges locally uniformly:
\[
 \phi_i\longrightarrow\phi\qquad\text{in }C_{\rm loc}(\R^k).
\]
For a unit vector $e$ and $h>0$, convexity gives
\[
 \frac{\phi_i(x)-\phi_i(x-he)}{h}
 \le\langle\nabla\phi_i(x),e\rangle
 \le\frac{\phi_i(x+he)-\phi_i(x)}{h}.
\]
Let $i\to\infty$ and then $h\rightarrow0$. Since $\phi$ is
smooth, this proves $\nabla\phi_i\to\nabla\phi$ uniformly on compact
sets. 
	
	\emph{Uniform local second derivative estimates.}
	The limit $\phi$ is smooth with positive definite Hessian. At a fixed $x_0$ let
	\[
	w_i(z)=\phi_i(z)-\phi_i(x_0)-\nabla\phi_i(x_0)\cdot(z-x_0),
	\qquad S_i(s)=\{w_i<s\}.
	\]
	Choose $r,\rho,h>0$ using the limit so that its analogous function is at least $4h$ on $\partial B(x_0,r)$ and at most $h/4$ on $B(x_0,\rho)$. Uniform convergence of potentials and gradients gives, for all sufficiently large $i$,
	\begin{equation}\label{eq:uniform-sections}
		B(x_0,\rho)\subset S_i(h/2)\subset S_i(2h)\Subset B(x_0,r).
	\end{equation}
	Here $A\Subset B$ means that the closure $\overline A$ is a compact subset of $B$.
	For the last inclusion, the value on $\partial B(x_0,r)$ is at least $3h$, and convexity rules out a point of $S_i(2h)$ beyond that sphere. The gradients on this ball lie in a common bounded target set. Hence
	\[
	\det \nabla^2\phi_i(z)=F_i(z):=\frac{e^{-\pi|z|^2}}{q_{E_i,t}(\nabla\phi_i(z))}
	\quad\hbox{and}\quad 0<\lambda\le F_i\le\Lambda,
	\]
   with constants independent of $i$.
	
	Use the interior estimates on nested sections, after affine normalization by John's ellipsoid. The inner and outer radii in \eqref{eq:uniform-sections} bound the normalizing linear maps and their inverses. The normalized $C^{1,\alpha}$ estimate gives a common gradient H\"older bound on a smaller section. Since $q_{E_i,t}$ has uniformly bounded first derivatives and a uniform positive lower bound on the target compact set, $F_i$ now has a common $C^\alpha$ bound there. The normalized $C^{2,\alpha}$ estimate then applies on a further smaller section. These are the section estimates stated in \cite[Theorems 2.16 and 2.20, with the localization following Theorem 2.22]{DPF}. In particular, they give uniform bounds near $x_0$. The lower determinant bound and the upper Hessian bound give a common positive lower Hessian bound as well. Thus the geometry of the sections and the constants in the local elliptic estimates are uniform along the parameter sequence.
	
	\emph{Higher derivatives and inverses.}
	Differentiating the logarithmic equation gives, for $v_i=\partial_\ell\phi_i$,
	\[
	\operatorname{tr}\!\left(
	(\nabla^2\phi_i)^{-1}\nabla^2v_i
	\right)
	+\nabla\log q_{E_i,t}(\nabla\phi_i)\cdot\nabla v_i
	=-2\pi z_\ell.
	\]
	The coefficients are uniformly elliptic and have common local H\"older bounds. Interior Schauder estimates and iteration give uniform bounds of every fixed order on smaller compact sets. Compactness and the already identified limit give convergence in every local $C^r$ norm. A finite covering argument gives the same conclusion on every compact subset of the source.
	
	For inverse maps, repeat the ball argument with source $q_{E_i,t}$ and Gaussian target. On each source compact set the lower bounds for $q_{E_i,t}$ give a common positive mass for the balls used above. Thus inverse images of compact sets remain in a common compact set. Uniform local convergence, the inverse function theorem, and differentiation of $U_i^{-1}\circ U_i=I$ now give convergence of all inverse derivatives. Local frame changes agree by uniqueness. 
\end{proof}

\section{Transportation proofs of the two inequalities}\label{sec:inequalities}

We prove both inequalities using the transport lemmas and the Grassmannian Ball--Barthe inequality. For the forward inequality, we first work on bounded balls and then use truncation and regularization. For the reverse inequality, we first consider uniformly tight Borel families of densities. We use this result again in Section~\ref{sec:reverse-equality}.

\subsection{The forward inequality}
First take the positive probability densities $a_E$ on $B_E(R)$ from Lemma~\ref{lem:ball-transport}, extended by zero outside. Put
\[
\Omega_R=\{x:\max_{E\in S}|P_Ex|<R\},\qquad
T(x)=\int_ST_{E,\rho}(P_Ex)\,d\nu.
\]
This set is open, convex and bounded, since
$\max_E|P_Ex|^2\ge(k/n)|x|^2$. Parameter continuity permits differentiation under the integral:
\[
DT(x)=\int_S\iota_E DT_{E,\rho}(P_Ex)P_E\,d\nu>0.
\]
Integration along line segments makes $T$ strictly monotone and injective.  The transport Jacobian identity gives
\[
 \log a_E(P_Ex)=-\log Z_\rho-\pi|T_{E,\rho}(P_Ex)|^2
                  +\log\det_E DT_{E,\rho}(P_Ex).
\]
Since $T(x)=\int_ST_{E,\rho}(P_Ex)\,d\nu(E)$, the energy identity and the determinant inequality give
\[
 \int_S|T_{E,\rho}(P_Ex)|^2\,d\nu(E)\ge|T(x)|^2,
 \qquad
 \exp\int_S\log\det_E DT_{E,\rho}(P_Ex)\,d\nu(E)\le\det DT(x).
\]
Using also $\nu(S)=n/k$, we obtain
\[
\exp\int_S\log a_E(P_Ex)\,d\nu
\le Z_\rho^{-n/k}e^{-\pi|T(x)|^2}\det DT(x).
\]
Consequently
\begin{equation}\label{eq:forward-ball}
	\int_{\Omega_R}\exp\int_S\log a_E(P_Ex)\,d\nu\,dx
	\le Z_\rho^{-n/k}\int_{T(\Omega_R)}e^{-\pi|z|^2}\,dz
	\le Z_\rho^{-n/k}.
\end{equation}
Outside $\overline{\Omega_R}$ some projection lies strictly outside its ball, and so do projections in a parameter neighborhood of positive measure. The product therefore vanishes there. The convex boundary has Lebesgue measure zero. Letting $\rho\to\infty$ proves the bound with constant one. The expression on the left is independent of $\rho$, so no convergence of the maps is needed.

Next suppose that a jointly continuous nonnegative family $f_E$
is supported in $\overline{B_E(R)}$. Choose a smooth radial probability
density $\eta$ supported in $B^k(1)$, and write
$\eta_\varepsilon(y)=\varepsilon^{-k}\eta(y/\varepsilon)$ in every
fiber. For $0<\varepsilon<1/2$ and $\delta>0$, set
\[
 f_{E,\varepsilon,\delta}(y)
 =\bigl((f_E*\eta_\varepsilon)(y)+\delta\bigr)
    \mathbf 1_{\overline{B_E(R+1)}}(y).
\]
Its restriction to the closed ball is positive and smooth in the
fiber. The radial mollifier is independent of the orthonormal frame,
and compactness gives the required parameter continuity. Normalize
each fiber, apply the ball result, and restore its mass:
\[
 \int_{\R^n}F_{f_{\varepsilon,\delta}}(x)\,dx
 \le\exp\int_S\log\!\left(\int_Ef_{E,\varepsilon,\delta}\right)
 \,d\nu(E).
\]
For fixed $\delta>0$, the mollified functions converge uniformly to
$f_E$ on the closed balls in local frames. Their logarithms converge
uniformly there because they are bounded below by $\delta$. The
ambient products have a common bounded support and are uniformly
bounded. The possible boundary of the ambient integration region has
Lebesgue measure zero. Letting $\varepsilon\rightarrow0$ therefore
gives the same inequality with
$f_{E,\delta}=(f_E+\delta)\mathbf1_{\overline{B_E(R+1)}}$.

As $\delta\rightarrow0$, the fiber values and masses decrease to
$f_E$ and $m_E$. Their logarithms are uniformly bounded above, so
monotone convergence applied after subtracting a common upper bound
also covers zero values and zero masses. The ambient products are
dominated on a common bounded set. We obtain the desired inequality
for the original compactly supported family.

For a general family choose radial cutoffs $\chi_R\rightarrow1^-$, equal
to one on $B(R)$ and zero outside $B(2R)$. Set
\[
 f_{E,R}=f_E\chi_R,\qquad m_{E,R}=\int_Ef_{E,R}.
\]
For $|x|\le R$, orthogonal projection gives $|P_Ex|\le R$. Hence
\[
 f_{E,R}(P_Ex)=f_E(P_Ex),\qquad F_{f_R}(x)=F_f(x).
\]
Since $m_{E,R}\le m_E$, the compact-support inequality gives
\[
\int_{B^n(R)}F_f(x)\,dx
\le\exp\int_S\log m_{E,R}\,d\nu\le\Mass(f).
\]
Let $R\to\infty$. This proves \eqref{eq:forward-main}, as well as the inequality under the weaker assumptions of Appendix~\ref{sec:extension}. 

\subsection{The reverse inequality for a uniformly tight family}
We prove an auxiliary assertion: the reverse bound holds for every density family from Section~\ref{sec:transport}, after multiplication by a positive mass field with integrable logarithm.

First let $a_E$ be a jointly Borel, fiber-continuous family of probability densities with common bounded support and a common upper bound $K$. Fiber continuity ensures the universal measurability of its envelope by Proposition~\ref{prop:measurability}. Use Lemma~\ref{lem:weak-transport}, and, on its common full-measure set of ambient points, set
\[
B(x)=\int_S\iota_E\nabla^2\phi_E(P_Ex)P_E\,d\nu,\qquad M(x)=\nabla^2_a\Phi(x).
\]
These matrices are finite, with $M(x)\succeq B(x)$. To see that
$B(x)$ is positive definite, suppose $\langle B(x)\xi,\xi\rangle=0$.
Then
\[
 0=\int_S\left\langle\nabla^2\phi_E(P_Ex)P_E\xi,P_E\xi\right\rangle
 \,d\nu(E).
\]
The integrand is nonnegative and the fiber Hessian is positive
definite for almost every $E$. Hence $P_E\xi=0$ almost everywhere.
By isotropy,
\[
 |\xi|^2=\int_S|P_E\xi|^2\,d\nu(E)=0.
\]
Thus $B(x)>0$, and $M(x)\succeq B(x)$ gives $M(x)>0$. The Jacobian equation gives
\[
\log\det_E\nabla^2\phi_E(P_Ex)\ge-\pi|x|^2-\log K.
\]
Its positive part is bounded by the trace, integrable in $E$ for almost every $x$ by \eqref{eq:weak-Hessian}. The displayed lower bound controls its negative part. Thus $\log\det \nabla^2\phi_E(P_Ex)$ is integrable in $E$. The Jacobian equation also makes $\log a_E(\nabla\phi_E(P_Ex))$ integrable there, so the logarithmic sums below contain no undefined cancellation. The measurable version of Lemma~\ref{lem:BB} applies. After altering the exceptional parameter values on a null set, the section $\nabla\phi_E(P_Ex)$ is bounded and has barycenter $T(x)$. Hence
\begin{align}
	\RB(a)(T(x))\det M(x)
	&\ge \exp\int_S\log a_E(\nabla\phi_E(P_Ex))\,d\nu\ \det B(x)\notag\\
	&\ge \exp\int_S\big(\log a_E(\nabla\phi_E(P_Ex))
	+\log\det_E\nabla^2\phi_E(P_Ex)\big)\,d\nu\notag\\
	&=e^{-\pi|x|^2}.\label{eq:reverse-weak-chain}
\end{align}
Lemma~\ref{lem:convex-area}, with the universally measurable function $H=\RB(a)$, gives $\int\RB(a)\ge1$. 

For a general density family $q$ satisfying the conditions above, put $r_{E,R}=\int q_E\chi_R$. Uniform tightness gives uniform convergence $r_{E,R}\to1$, so $r_{E,R}\ge1/2$ for all large $R$. The normalized compactly supported density $a_{E,R}=q_E\chi_R/r_{E,R}$ has a common upper bound. Monotonicity of the envelope yields
\[
\int\RB(q)\ge\int\RB(q\chi_R)
\ge\exp\int_S\log r_{E,R}\,d\nu\longrightarrow1.
\]
 Rescaling proves the auxiliary assertion and \eqref{eq:reverse-main}. Since $\RBone\ge\RB$, its version follows as well.

For $f_E(y)=e^{-\pi|y|^2}$, the vector identity gives
$\RB(f)(x)=e^{-\pi|x|^2}$, attained by $v(E)=P_Ex$. The same identity holds for $\RBone$. Both constants are optimal.

\subsection{A smooth reverse chain}
Let $q_{E,t}$ be a heat-regularized admissible probability family and let $U_E$ be the transports of Lemma~\ref{lem:heat-transport}. Define
\[
T(x)=\int_SU_E(P_Ex)\,d\nu,\qquad
M(x)=DT(x)=\int_S\iota_EDU_E(P_Ex)P_E\,d\nu.
\]
Each derivative $\nabla U_E$ is positive definite, and isotropy
therefore gives $M(x)>0$. If $x\ne y$, integration along the segment
from $y$ to $x$ yields
\[
 \langle T(x)-T(y),x-y\rangle
 =\int_0^1\langle M(y+s(x-y))(x-y),x-y\rangle\,ds>0.
\]
Thus $T$ is injective. The parameter estimates also show that $T$ is
smooth and each section $E\mapsto U_E(P_Ex)$ is continuous and
bounded. Set
\[
{\mathcal D}(x)=\log\det M(x)-\int_S\log\det_EDU_E(P_Ex)\,d\nu\ge0.
\]
Change of variables and the transport Jacobian equation imply
\begin{equation}\label{eq:smooth-reverse-chain}
	\int\RB(q_t)\ge
	\int_{\R^n}\RB(q_t)(T(x))\det M(x)\,dx
	\ge\int_{\R^n}e^{-\pi|x|^2}e^{{\mathcal D}(x)}\,dx\ge1.
\end{equation}
If $\int\RB(q_t)=1$, the continuous nonnegative function $\mathcal D$ vanishes everywhere. Lemma~\ref{lem:BB} therefore gives
\begin{equation}\label{eq:reverse-map-equality}
	M(x)P_E=P_EM(x),\qquad DU_E(P_Ex)=M(x)|_E
	\quad(x\in\R^n,\ E\in S).
\end{equation}

\section{Equality in the forward inequality}\label{sec:forward-equality}

By homogeneity, we may normalize each fiber so that $\int_Eq_E=1$. We first show that Gaussian convolution preserves extremizers. A differential identity then determines the structure of the smooth extremizers, and a fixed-time argument recovers the original family. Put
\[
F_0(x)=\exp\int_S\log q_E(P_Ex)\,d\nu,\qquad
G_t(x)=\exp\int_S\log q_{E,t}(P_Ex)\,d\nu .
\]
The heat-regularized family is admissible, so $\int G_t\le1$.

\begin{lemma}[Exact forward smoothing]\label{lem:forward-smoothing}
	For every $t>0$,
	\begin{equation}\label{eq:forward-convolution}
		G_t\ge F_0*\gamma_{n,t}.
	\end{equation}
	If $\int F_0=1$, then
	\begin{equation}\label{eq:forward-heat-identity}
		G_t=F_0*\gamma_{n,t}\quad\text{everywhere},\qquad\int G_t=1.
	\end{equation}
\end{lemma}
\begin{proof}
	Fix $t>0$ and $x\in\R^n$. Apply the forward inequality to
\[
 y\longmapsto q_E(y)\gamma_{k,t}(P_Ex-y),\qquad y\in E.
\]
Its mass is
\[
 \int_Eq_E(y)\gamma_{k,t}(P_Ex-y)\,dy
 =q_{E,t}(P_Ex)>0.
\]
The mass is continuous in $E$ by $L^1$ parameter continuity, so this
family is admissible. Isotropy gives
	\[
	\exp\int_S\log\gamma_{k,t}(P_E(x-z))\,d\nu
	=\gamma_{n,t}(x-z).
	\]
	The forward inequality for this family is
	\[
	\int_{\R^n}F_0(z)\gamma_{n,t}(x-z)\,dz
	\le\exp\int_S\log q_{E,t}(P_Ex)\,d\nu(E)=G_t(x).
	\]
	This proves \eqref{eq:forward-convolution}. Under equality at time zero,
	\[
	1=\int F_0*\gamma_{n,t}\le\int G_t\le1.
	\]
	The first and last terms are equal, so the nonnegative function
$G_t-F_0*\gamma_{n,t}$ has integral zero. It vanishes almost everywhere.
Both terms are continuous. $G_t$ is continuous by the parameter
estimates, and convolution of the integrable function $F_0$ with the
Gaussian is continuous. Their equality therefore holds everywhere.
\end{proof}

For $t>0$, $E\in S$, and $y\in E$, write $L_t=\log G_t$ and $s_E(t,y)=\nabla_E\log q_{E,t}(y)\in E$. We regard $s_E$ as an ambient vector in $\R^n$. On a compact time interval away from zero and a bounded spatial set, the heat-regularized densities $q_{E,t}$ have uniformly bounded derivatives and a common positive lower bound. The derivative bounds follow by convolution with probability densities. The lower bound follows by retaining a common ball carrying at least half the mass. Thus the logarithms and all their required derivatives are uniformly bounded there, and differentiation commutes with the finite parameter integral. Since each $q_{E,t}$ solves the heat equation, the logarithmic
heat equation on $E$ gives
\[
 \partial_t\log q_{E,t}
 =\Delta_E\log q_{E,t}+|s_E|^2.
\]
Moreover,
\[
 \nabla L_t(x)=\int_Ss_E(t,P_Ex)\,d\nu(E),\qquad
 \Delta L_t(x)=\int_S\Delta_E\log q_{E,t}(P_Ex)\,d\nu(E).
\]
Since $G_t=e^{L_t}$, the chain rule and Lemma~\ref{lem:energy} give
\begin{align}
	(\partial_t-\Delta)G_t
	&=G_t\left(\int_S|s_E(t,P_Ex)|^2\,d\nu
	-\left|\int_Ss_E(t,P_Ex)\,d\nu\right|^2\right)\notag\\
	&=G_t\int_S|s_E(t,P_Ex)-P_E\nabla L_t(x)|^2\,d\nu.\label{eq:forward-defect}
\end{align}
If $q$ is an extremizer, \eqref{eq:forward-heat-identity} makes the left side zero. Joint continuity and the support condition imply
\begin{equation}\label{eq:score-compatibility}
	\nabla_E\log q_{E,t}(P_Ex)=P_E\nabla L_t(x)
	\quad(E\in S,\ x\in\R^n).
\end{equation}
Differentiate \eqref{eq:score-compatibility} with respect to $x$. Since the Hessian on $E$ is symmetric,
\[
 \nabla_E^2\log q_{E,t}(P_Ex)P_E=P_E\nabla^2L_t(x),
 \qquad P_E\nabla^2L_t(x)=\nabla^2L_t(x)P_E.
\]

Fix $t=t_*>0$. Lemma~\ref{lem:rigidity} applies to $L_{t_*}$. Integrability of $G_{t_*}$ forces its constant quadratic form on $D$ to be negative definite: a positive eigenvalue, or a zero eigenvalue with any linear term, makes the exponential nonintegrable along that direction. Hence there are a positive definite covariance $\Sigma_*$ on $D$, a vector $a\in D$, and probability densities $h_{j,*}$ such that
\begin{equation}\label{eq:forward-global-factor}
	G_{t_*}(x)=\varphi_{D,a,\Sigma_*}(x_D)\prod_jh_{j,*}(x_j).
\end{equation}
Here $\varphi_{L,c,\Sigma}$ is the Gaussian density defined in Lemma~\ref{lem:deconvolution}. The operator $\Sigma_*$ commutes with every $P_{E_D}$. By \eqref{eq:score-compatibility} and connectedness of each fiber, normalization gives
\begin{equation}\label{eq:forward-fiber-factor}
	q_{E,t_*}(y)=
	\varphi_{E_D,P_{E_D}a,\Sigma_*|_{E_D}}(P_Dy)
	\prod_{F_j\subset E}h_{j,*}(P_{F_j}y).
\end{equation}

\subsection{Recovery at a fixed positive time}
Fix the time $t_*>0$. Since $p_0=F_0$ is a probability density and $G_{t_*}=p_0*\gamma_{n,t_*}$, Lemma~\ref{lem:deconvolution} applies. To identify the factors on different fibers, let $p_{0,D}$ and $h_j$ be the marginals of $p_0$ on $D$ and $F_j$. By Fubini, they are probability densities. With
\[
 \widehat p(\xi)=\int_{\R^n}e^{i\langle\xi,x\rangle}p(x)\,dx,
\]
the Fourier identity on $D$ gives
\[
\widehat p_{0,D}(\xi_D)
=\exp\left(i\langle a,\xi_D\rangle
-\tfrac12\langle(\Sigma_*-2t_*I_D)\xi_D,\xi_D\rangle\right).
\]
Set $\Sigma_0=\Sigma_*-2t_*I_D$. Boundedness of a characteristic function forces $\Sigma_0\succeq0$. A nonzero kernel vector would make $p_{0,D}$ supported on a proper affine hyperplane, contradicting absolute continuity. Thus $\Sigma_0>0$ on $D$.

Similarly,
$\widehat h_{j,*}(\xi_j)=\widehat h_j(\xi_j)e^{-t_*|\xi_j|^2}$.
For arbitrary $\xi$, cancellation of the Gaussian factors gives
\[
\widehat p_0(\xi)=
e^{i\langle a,\xi_D\rangle-\langle\Sigma_0\xi_D,\xi_D\rangle/2}
\prod_j\widehat h_j(\xi_j).
\]
Fix $j$ and $E\in S_j$, and let the probability density of the
$F_j$ marginal of $q_E$ be denoted by $p_{E,j}$. The factorization
\eqref{eq:forward-fiber-factor} gives
\[
 \widehat p_{E,j}(\xi_j)e^{-t_*|\xi_j|^2}
 =\widehat h_{j,*}(\xi_j)
 =\widehat h_j(\xi_j)e^{-t_*|\xi_j|^2}.
\]
The Gaussian multiplier is never zero. Thus
\[
 p_{E,j}=h_j\quad\text{almost everywhere on }F_j.
\]
This identifies the same marginal in every fiber containing $F_j$.
Applying the same cancellation to all coordinates in
\eqref{eq:forward-fiber-factor}, and using uniqueness of characteristic
functions, gives the joint factorization below. Consequently
\begin{equation}\label{eq:forward-recovered}
	q_E(y)=\varphi_{E_D,P_{E_D}a,\Sigma_0|_{E_D}}(P_Dy)
	\prod_{F_j\subset E}h_j(P_{F_j}y)
	\quad\text{almost everywhere in }E.
\end{equation}
Each recovered factor is an actual marginal of $p_0$. In particular, its positivity is obtained from the marginal, rather than from a formal inverse Fourier transform.

The continuous-factor observation following Lemma~\ref{lem:rigidity} gives continuous versions of all $h_j$. The identity then holds pointwise in every fiber. Put
\[
A=(2\pi)^{-1}\Sigma_0^{-1},\qquad
\theta_E=m_E\det(A|_{E_D})^{1/2}.
\]
The restrictions vary continuously (their ranks are locally constant on $S$), so $\theta_E$ is positive and continuous. This gives \eqref{eq:forward-form}.

\subsection{Sufficiency}
Assume \eqref{eq:forward-form}, with $\int h_j=1$. From
\eqref{eq:block-identities},
\[
F_f(x)=
e^{\int\log\theta_E\,d\nu}
e^{-\pi\langle A(x_D-a),x_D-a\rangle}\prod_jh_j(x_j).
\]
Also,
\[
m_E=\theta_E\det(A|_{E_D})^{-1/2},\qquad
\int_S\log\det(A|_{E_D})\,d\nu=\tr_D\log A=\log\det_DA.
\]
Since the $h_j$ have integral one, the first identity gives
\[
 \BL(f)=e^{\int_S\log\theta_E\,d\nu(E)}(\det_DA)^{-1/2}.
\]
The identities for $m_E$ and $\det(A|_{E_D})$ give the same value for $\Mass(f)$. These formulas also hold when the factors vanish, with $\log0=-\infty$. This proves Theorem~\ref{thm:forward}.

\section{Equality in the reverse inequality}\label{sec:reverse-equality}

By homogeneity, we again normalize $\int_Eq_E=1$. We first establish a convolution comparison showing that Gaussian regularization preserves reverse extremizers. A finite-family version of this method is used in \cite{BKX}. Here the auxiliary inequality for Borel families allows the sections in the supremum to depend measurably on the parameter.

\begin{lemma}[Exact reverse smoothing]\label{lem:reverse-smoothing}
	For every $t>0$,
	\begin{equation}\label{eq:reverse-convolution}
		\RB(q_t)(x)\le(\RB(q)*\gamma_{n,t})(x).
	\end{equation}
	If $\int\RB(q)=1$, then $\int\RB(q_t)=1$ for every $t>0$, and the two sides of \eqref{eq:reverse-convolution} agree almost everywhere.
\end{lemma}
\begin{proof}
	Fix a bounded section $v$ with $\int v\,d\nu=x$, choose an everywhere bounded Borel representative as in Proposition~\ref{prop:measurability}, and set
	\[
	b_E^v(y)=q_E(y)\gamma_{k,t}(v(E)-y).
	\]
	Although $v$ need not be continuous, this is a Borel family, continuous in each fiber. Write $|v(E)|\le N$ and choose $R_0$ so that every $q_E$ gives $B_E(R_0)$ mass at least $1/2$. Its mass $\beta_E^v=q_{E,t}(v(E))$ satisfies
	\[
	\tfrac12(4\pi t)^{-k/2}e^{-(N+R_0)^2/(4t)}
	\le\beta_E^v\le(4\pi t)^{-k/2}.
	\]
	Thus $\log\beta_E^v$ is integrable. Also
	$b_E^v/\beta_E^v\le C_{N,t,R_0}q_E$, which proves uniform tightness and local uniform bounds for this normalized family. The auxiliary reverse inequality, already proved in Section~\ref{sec:inequalities} without the equality classification or the smooth transport lemma, therefore gives
	\[
	\exp\int_S\log q_{E,t}(v(E))\,d\nu
	\le\int_{\R^n}\RB(b^v)(z)\,dz.
	\]
	For any bounded section $w$ with barycenter $z$, Lemma~\ref{lem:energy} yields
	\[
	\exp\int_S\log\gamma_{k,t}(v(E)-w(E))\,d\nu
	\le\gamma_{n,t}(x-z).
	\]
	Taking the supremum over such $w$ gives
	\[
	\RB(b^v)(z)\le\RB(q)(z)\gamma_{n,t}(x-z).
	\]
	Hence
	\[
	\exp\int_S\log q_{E,t}(v(E))\,d\nu(E)
	\le\int_{\R^n}\RB(q)(z)\gamma_{n,t}(x-z)\,dz.
	\]
	Taking the supremum over $v$ proves \eqref{eq:reverse-convolution}. If $\int\RB(q)=1$, the reverse inequality gives
	\[
	1\le\int\RB(q_t)\le\int\RB(q)*\gamma_{n,t}=1.
	\]
	This proves exact extremality at every positive time. No continuity of the reverse envelope is needed to deduce its almost-everywhere equality with the convolution.
\end{proof}

\subsection{Structure of a smooth reverse extremizer}
Fix $t_*>0$ and suppose that $\int\RB(q)=1$. By the preceding lemma, $q_{t_*}$ is a reverse extremizer. Use the smooth maps $U_E$ from Gaussian measure to $q_{E,t_*}$ and the notation of \eqref{eq:smooth-reverse-chain}. Equality yields \eqref{eq:reverse-map-equality} everywhere. Lemma~\ref{lem:rigidity}, applied to the potential of $T=\int U_E\circ P_E\,d\nu$, gives
\[
T(x)=Bx_D+c+\sum_jt_j(x_j),
\]
where $B>0$ on $D$, $B$ commutes with all $P_{E_D}$, and each $t_j$ is a convex-gradient map with positive definite derivative. Absorb $t_j(0)$ into the constants so that $t_j(0)=0$. Since the derivatives of $U_E$ are the corresponding restrictions of $DT$, integration on each fiber gives
\begin{equation}\label{eq:reverse-block-map}
	U_E(y)=BP_Dy+\sum_{F_j\subset E}t_j(P_{F_j}y)+b_E,\qquad b_E=U_E(0).
\end{equation}
The section $b_E$ is continuous and bounded on $S$. Each $U_E$ is onto. Every $F_j$ occurs in a fiber, so each $t_j$ is onto as well; strict monotonicity makes it a diffeomorphism.

Write $y=P_Dy+\sum_{F_j\subset E}P_{F_j}y$. By orthogonality,
\[
 e^{-\pi|y|^2}=e^{-\pi|P_Dy|^2}
 \prod_{F_j\subset E}e^{-\pi|P_{F_j}y|^2}.
\]
Applying the block map \eqref{eq:reverse-block-map} to these independent Gaussian factors gives
\begin{equation}\label{eq:reverse-smooth-factor}
	q_{E,t_*}(y)=
	\varphi_{E_D,P_Db_E,\Sigma_*|_{E_D}}(P_Dy)
	\prod_{F_j\subset E}h_{j,*}(P_{F_j}(y-b_E)),
	\qquad \Sigma_*=(2\pi)^{-1}B^2,
\end{equation}
where $h_{j,*}$ is the probability density obtained by pushing
$e^{-\pi|z|^2}\,dz$ on $F_j$ through $t_j$. These factors are common to all fibers containing $F_j$.

\subsection{Deconvolution and parameter consistency}
Set $\Sigma_0=\Sigma_*-2t_*I_D$. Apply Lemma~\ref{lem:deconvolution} in each fiber, after translating by $-b_E$. This is legitimate for every $E$, since smooth extremality and the transport identities hold on all of $S$. In particular, restricting the characteristic-function identity
$\widehat q_{E,t_*}=\widehat q_Ee^{-t_*|\xi|^2}$ to $E_D$ shows that the original $E_D$ marginal has characteristic function
\[
\exp\big(i\langle P_Db_E,\xi_D\rangle
-\tfrac12\langle\Sigma_0\xi_D,\xi_D\rangle\big).
\]
Boundedness of its characteristic function and absolute continuity of the marginal imply
$\Sigma_0|_{E_D}>0$ whenever $E_D\ne0$. This forces $\Sigma_0>0$ on all of $D$. Indeed, $\Sigma_0$ commutes with every $P_{E_D}$; the projection of an eigenvector with nonpositive eigenvalue would either vanish or be a nonpositive-eigenvalue vector in $E_D$. The restrictions forbid the latter, whereas isotropy forbids all projections vanishing.

For each $j$, choose one reference fiber $E(j)\in S_j$. Let $h_j$ be the $F_j$ marginal of $q_{E(j)}$, translated by $-P_{F_j}b_{E(j)}$. Explicitly, if $r_{E,j}$ denotes that marginal, then $h_j(z)=r_{E(j),j}(z+P_{F_j}b_{E(j)})$ almost everywhere. It is an actual probability density, and
\[
\widehat h_{j,*}(\xi_j)=\widehat h_j(\xi_j)e^{-t_*|\xi_j|^2}.
\]
Substitute this identity and the covariance relation into
\eqref{eq:reverse-smooth-factor}. Since the fiber decomposition is orthogonal, all heat factors cancel, and characteristic-function uniqueness gives
\begin{equation}\label{eq:reverse-recovered}
	q_E(y)=
	\varphi_{E_D,P_Db_E,\Sigma_0|_{E_D}}(P_Dy)
	\prod_{F_j\subset E}h_j(P_{F_j}(y-b_E))
	\quad\text{almost everywhere in }E.
\end{equation}
The reference fiber determines the same factor for every fiber containing $F_j$. The translations $b_E$ come from the single fixed time $t_*$ and remain continuous and bounded. No first-moment assumption on the independent factors is needed.

For $E\in S_j$, write
\[
 \widetilde r_{E,j}(z)=r_{E,j}(z+P_{F_j}b_E),\qquad z\in F_j.
\]
The smooth factorization gives
\[
 \widehat{\widetilde r}_{E,j}(\xi_j)e^{-t_*|\xi_j|^2}
 =\widehat h_{j,*}(\xi_j)
 =\widehat h_j(\xi_j)e^{-t_*|\xi_j|^2}.
\]
The Gaussian multiplier is nonzero, so uniqueness of characteristic
functions gives $\widetilde r_{E,j}=h_j$ almost everywhere. This holds
for every $E\in S_j$, including parameters with zero individual
$\nu$-mass. Choose continuous versions of the factors using the fixed
reference fibers. The fiber identities then extend from almost
everywhere to everywhere by continuity.

By the continuous-slice argument of Section 2, each $h_j$ has a continuous version. The resulting product is continuous in each fiber, so \eqref{eq:reverse-recovered} holds pointwise. Taking
\[
A=(2\pi)^{-1}\Sigma_0^{-1},\qquad
\theta_E=m_E\det(A|_{E_D})^{1/2}
\]
gives \eqref{eq:reverse-form}. It remains to determine exactly which independent factors are allowed.

\subsection{Log-concavity of the shared factors}
Let $\bar b=\int_S b_E\,d\nu$. The bounded section
$v(E)=b_E+P_E(x-\bar b)$ is feasible in \eqref{eq:reverse-envelope}.
Formula \eqref{eq:reverse-recovered} therefore implies
\begin{equation}\label{eq:reverse-baseline}
	\RB(q)(x)\ge
	H(x):=\varphi_{D,P_D\bar b,\Sigma_0}(x_D)
	\prod_jh_j(x_j-P_{F_j}\bar b).
\end{equation}
The Gaussian normalization follows from
\[
 \int_SP_{E_D}\,d\nu(E)=I_D,\qquad
 \Sigma_0P_{E_D}=P_{E_D}\Sigma_0.
\]
Each factor has weight $\nu(S_j)=1$. Since all factors in $H$ are probability densities,
\[
 \int_{\R^n}H(x)\,dx=1.
\]

Suppose $F_j$ is shared. Choose a Borel set $A_0\subset S_j$ with
$\lambda=\nu(A_0)\in(0,1)$. We claim
\begin{equation}\label{eq:fixed-lambda-concavity}
	h_j(\lambda u+(1-\lambda)v)
	\ge h_j(u)^\lambda h_j(v)^{1-\lambda}
	\quad(u,v\in F_j).
\end{equation}
If this fails, the right side is positive. Put $r=\lambda u+(1-\lambda)v$.
In the section used for \eqref{eq:reverse-baseline}, replace the $F_j$ component relative to $b_E$ by $u+z$ on $A_0$ and by $v+z$ on $S_j\setminus A_0$. Leave the other components unchanged. The new $F_j$ barycenter and factor are
\[
 P_{F_j}\bar b+r+z,
 \qquad h_j(u+z)^\lambda h_j(v+z)^{1-\lambda}.
\]
Write $x_{\ne j}=(x_D,(x_i)_{i\ne j})$, and let $K$ be the product of the Gaussian and all the factors of $H$ other than $h_j$. Using both this section and the original section, we obtain
\[
 \RB(q)(x)\ge K(x_{\ne j})\max\left\{
 h_j(r+z),\ h_j(u+z)^\lambda h_j(v+z)^{1-\lambda}\right\},
 \qquad x_j=P_{F_j}\bar b+r+z.
\]
The second term is strictly larger near $z=0$, by continuity. As $\int K=\int h_j=1$, integration gives $\int\RB(q)>1$. This contradiction proves \eqref{eq:fixed-lambda-concavity}.

We now pass from the fixed weight $\lambda$ to every weight
$\theta\in[0,1]$. Iterate \eqref{eq:fixed-lambda-concavity} along a
binary tree of depth $m$. Its positive leaf weights sum to one and
each is at most $\max(\lambda,1-\lambda)^m$. Assigning $u$ or $v$
to the leaves gives
\[
 h_j(\theta_m u+(1-\theta_m)v)
 \ge h_j(u)^{\theta_m}h_j(v)^{1-\theta_m},
\]
where $\theta_m$ is the sum of the weights assigned to $u$. For any
$\theta\in[0,1]$, select leaves in order until their total first
exceeds $\theta$. The resulting $\theta_m$ satisfies
\[
 |\theta_m-\theta|\le\max(\lambda,1-\lambda)^m\longrightarrow0.
\]
If $h_j(u),h_j(v)>0$, continuity gives the desired inequality for
$\theta$. If either endpoint value is zero, the inequality is
immediate for $0<\theta<1$. The endpoint cases are identities. Thus
$h_j$ is log-concave.

If $\nu|_{S_j}$ is a point mass, write it as $\delta_{E_0}$;
this atom has mass one because $\nu(S_j)=1$. For $u\in E_0$,
isotropy gives
\[
 |u|^2=\int_S|P_Eu|^2\,d\nu(E)
      =|u|^2+\int_{S\setminus\{E_0\}}|P_Eu|^2\,d\nu(E).
\]
Thus $P_Eu=0$ almost everywhere on $S\setminus\{E_0\}$. Continuity
and the support condition extend this relation to every $E$ in the
support of the remaining measure. Also, $S_j$ is relatively open and
$\nu|_{S_j}=\delta_{E_0}$, so $S_j=\{E_0\}$. It follows that all
other supporting subspaces are orthogonal to $E_0$, and the complete
common eigenspace with eigenvalue one on $E_0$ is
\[
 F_j=E_0.
\]
The parameter measure on this block is a single unit atom, so no
log-concavity condition is needed for its factor.

\subsection{Sufficiency and integrable sections}
Assume \eqref{eq:reverse-form}. Let $\bar b=\int b_E\,d\nu$, and let $v$ be any admissible section with barycenter $x$. First we check the separate positive logarithmic parts. For a shared block, the continuous log-concave probability density $h_j$ has a nonempty open positivity set. A supporting affine function for the concave function $\log h_j$ at an interior point is a global upper bound. Evaluated on an $L^1$ section, its positive part is integrable. For an unshared block the parameter measure is a unit atom, so its contribution has finite positive part as well. Finally, $\log\theta_E$ is bounded. Thus none of the following decompositions creates an $\infty-\infty$ expression.

On $D$, if the quadratic energy is finite, the vector identity after application of $A^{1/2}$ gives
\[
\int_S\langle A(P_Dv(E)-P_Db_E),P_Dv(E)-P_Db_E\rangle\,d\nu
\ge\langle A(x_D-P_D\bar b),x_D-P_D\bar b\rangle .
\]
If the left side is infinite, the preceding positive-part bounds make the full product zero, so the resulting product bound is automatic.
For each $j$, the section condition implies
\[
 P_{F_j}v(E)=P_{F_j}b_E=0\qquad(E\notin S_j)
\]
outside a parameter null set. Hence
\[
 \int_{S_j}P_{F_j}(v(E)-b_E)\,d\nu(E)
 =x_j-P_{F_j}\bar b.
\]
For a shared $F_j$, Jensen's inequality for the concave extended
function $\log h_j$ and the probability measure $\nu|_{S_j}$ gives
\[
\int_{S_j}\log h_j(P_{F_j}(v(E)-b_E))\,d\nu
\le\log h_j(x_j-P_{F_j}\bar b).
\]
This remains valid for integrable sections: $\log h_j$ is upper semicontinuous and concave as an extended function, and the affine upper bound just verified permits the extended Jensen inequality, with $-\infty$ allowed. For an unshared $F_j$, its unit atom gives the same relation as an identity.

It follows that both envelopes are bounded above by
\begin{equation}\label{eq:reverse-envelope-extremizer}
	e^{\int\log\theta_E\,d\nu}
	e^{-\pi\langle A(x_D-P_D\bar b),x_D-P_D\bar b\rangle}
	\prod_jh_j(x_j-P_{F_j}\bar b).
\end{equation}
The bounded section $v(E)=b_E+P_E(x-\bar b)$ has barycenter $x$ and attains this value. Thus both envelopes equal \eqref{eq:reverse-envelope-extremizer}. Integration gives
\[
 \int_{\R^n}\RB(f)=\int_{\R^n}\RBone(f)
 =e^{\int_S\log\theta_E\,d\nu(E)}(\det_DA)^{-1/2}
 =\Mass(f),
\]
where the last equality uses the determinant identity from the forward proof.

Conversely, if $\int\RBone(q)=1$, then
$1\le\int\RB(q)\le\int\RBone(q)=1$.
The preceding necessity proof applies to $\RB$, yielding exactly the same structure. This completes the proof of Theorem~\ref{thm:reverse}.

\section{Consequences and comparisons}\label{sec:consequences}

We now apply the equality conditions to several special cases. We first consider full support and irreducible families of projections, and then characterize equality for a common radial profile. Finally, we use homogeneous profiles to obtain volume inequalities and their equality conditions.

\begin{corollary}[Full support]\label{cor:full}
	Suppose $1\le k<n$ and $S=\G$. Forward equality holds precisely for
	\[
	f_E(y)=\theta_Ee^{-\pi\lambda|y-P_Ea|^2},
	\qquad \lambda>0,\ a\in\R^n.
	\]
	Reverse equality holds precisely for
	\[
	f_E(y)=\theta_Ee^{-\pi\lambda|y-b_E|^2},
	\qquad \lambda>0,
	\]
	where $b_E\in E$ is a continuous section. The amplitudes are positive and continuous. The same conclusions hold if the supporting projections have no nontrivial common invariant subspace.
\end{corollary}
\begin{proof}
	Suppose first that the supporting projections have no nontrivial common invariant subspace. Every common $0/1$ eigenspace is invariant, so it is either $\{0\}$ or $\R^n$. The latter is impossible because $1\le k<n$. Hence
	\[
	 F_1\oplus\cdots\oplus F_s=\{0\},\qquad D=\R^n.
	\]
	If a symmetric operator $A$ commutes with all the projections, each eigenspace of $A$ is invariant under them. Therefore $A=\lambda I_n$. The main theorems give the stated forms.

	For $S=G(n,k)$, let $V$ be a nonzero subspace invariant under all
$P_E$, and choose $0\ne v\in V$. For any $w$ with
$\langle v,w\rangle\ne0$, choose a $(k-1)$-dimensional subspace of
$\operatorname{span}\{v,w\}^{\perp}$ when $v,w$ are independent; if
they are dependent, choose it in $w^\perp$. Together with $w$ it spans
a $k$-plane $E$ satisfying
\[
 P_Ev=\frac{\langle v,w\rangle}{|w|^2}w.
\]
Since $V$ is invariant and the scalar is nonzero, $w\in V$. Such
vectors $w$ span $\R^n$, so $V=\R^n$.
\end{proof}

\begin{example}[Different translation conditions]
	In $\R^2$, take unit vectors $u_1,u_2,u_3$ at angles $0,2\pi/3,4\pi/3$, put $E_i=\R u_i$, and give each line mass $2/3$. For
	$q_i(su_i)=e^{-\pi(s-\mu_i)^2}$, the reverse inequality is an equality for every $(\mu_1,\mu_2,\mu_3)$, with envelope
	\[
	\RB(q)(x)=e^{-\pi|x-\bar b|^2},\qquad
	\bar b=\tfrac23\sum_i\mu_i u_i .
	\]
	The forward integral is
	\[
	\BL(q)=\exp\left[-\pi\left(\tfrac23\sum_i\mu_i^2-|\bar b|^2\right)\right].
	\]
	For $(\mu_1,\mu_2,\mu_3)=(0,0,1)$ it is $e^{-2\pi/9}<1$.
\end{example}

\begin{example}[Shared independent factors]
	In $\R^3$, let
\[
 E_{12}=\operatorname{span}\{e_1,e_2\},\quad
 E_{13}=\operatorname{span}\{e_1,e_3\},\quad
 E_{23}=\operatorname{span}\{e_2,e_3\},\quad
 \nu=\tfrac12(\delta_{E_{12}}+\delta_{E_{13}}+\delta_{E_{23}}).
\]
Then $F_j=\R e_j$ for $j=1,2,3$, $D=\{0\}$, and every $F_j$ is shared. Let $h_j$ be a continuous probability density on $\R$. With $q_{ij}(x_ie_i+x_je_j)=h_i(x_i)h_j(x_j)$, we have
	\[
	q_{12}=h_1h_2,\qquad q_{13}=h_1h_3,\qquad q_{23}=h_2h_3
	\]
	This family is a forward extremizer for arbitrary continuous probability densities $h_j$. It is a reverse extremizer if and only if every $h_j$ is log-concave. Thus a non-log-concave two-bump density in one factor separates the two equality classes.
\end{example}

\begin{corollary}[A common radial profile]\label{cor:radial}
	Let $q_E(y)=\rho(|y|)$ be the same nonnegative continuous radial probability density on every fiber. Equality in either main inequality holds if and only if at least one of the following holds:
	\begin{enumerate}
		\item $\rho(r)=\lambda^{k/2}e^{-\pi\lambda r^2}$ for some $\lambda>0$;
		\item $k$ divides $n$ and $\nu=\sum_{i=1}^{n/k}\delta_{E_i}$, where the $E_i$ are mutually orthogonal $k$-planes.
	\end{enumerate}
\end{corollary}
\begin{proof}
	Suppose a radial probability density factors over two nonzero orthogonal subspaces. Its characteristic function has the form $\Psi(|\xi|^2)$. Independence gives
	\[
	 \Psi(s+t)=\Psi(s)\Psi(t),\qquad
	 \Psi(0)=1,\qquad \Psi(s)=\Psi(s/2)^2
	 \quad(s,t\ge0).
	\]
	By continuity and $\Psi(0)=1$, the function is positive near
zero. Since $\Psi(s)=\Psi(s/2^m)^{2^m}$, it is positive for every
$s\ge0$. The function $-\log\Psi$ is therefore continuous and
additive on $[0,\infty)$, so $\Psi(s)=e^{-cs}$ for some real $c$.
Boundedness of a characteristic function gives $c\ge0$. If $c=0$,
uniqueness of characteristic functions gives a point mass at the
origin, contrary to absolute continuity. Thus $c>0$ and the density
is Gaussian.
	
	Apply either equality classification. A non-Gaussian radial fiber cannot have a nonzero Gaussian block and another factor, cannot be entirely Gaussian, and cannot split into two independent factors. It must therefore be a single complete independent space $F_j$ of dimension $k$. The support consists of mutually orthogonal such spaces, and isotropy forces them to span $\R^n$ with unit atomic weights. Conversely, the Gaussian case follows from isotropy and the orthogonal case from Fubini.
\end{proof}

\subsection{Homogeneous functionals and volume bounds}
For $p\ge1$ and $x\in\R^n$, define
\[
H_p(x)=\int_S|P_Ex|^p\,d\nu,\qquad
I_p(x)=\inf_{\substack{v\in L^\infty,\ v(E)\in E\\\int v\,d\nu=x}}
\int_S|v(E)|^p\,d\nu .
\]
Both are finite, continuous, positive away from zero, convex and homogeneous of degree $p$. For $I_p$, finiteness follows by $v(E)=P_Ex$, convexity by interpolation of sections, and coercivity by
\[
I_p(x)\ge \nu(S)^{1-p}|x|^p
\]
(with the usual interpretation at $p=1$). Put
\[
K_p^+=\{H_p\le1\},\qquad K_p^-=\{I_p\le1\},\qquad
\alpha_{n,k,p}=
\frac{[\omega_k\Gamma(1+k/p)]^{n/k}}{\Gamma(1+n/p)},
\]
where $\omega_d$ is the volume of the Euclidean unit ball in $\R^d$.

\begin{corollary}\label{cor:volume}
	For $p\ge1$,
	\begin{equation}\label{eq:volumes}
		|K_p^+|\le\alpha_{n,k,p}\le|K_p^-|.
	\end{equation}
	If $p=2$, both bodies are the Euclidean unit ball and both bounds are equalities. If $p\ne2$, equality in either bound holds precisely for the orthogonal unit-atom measures in Corollary~\ref{cor:radial}. In particular, both bounds are strict when $p\ne2$ and $k$ does not divide $n$.
\end{corollary}
\begin{proof}
	Take $f_E(y)=e^{-|y|^p}$. Polar coordinates give
	\[
	 \int_Ef_E(y)\,dy
	 =k\omega_k\int_0^\infty e^{-r^p}r^{k-1}\,dr
	 =\omega_k\Gamma(1+k/p).
	\]
	By the definitions of $H_p$ and $I_p$,
	\[
	 F_f(x)=e^{-H_p(x)},\qquad
	 \RB(f)(x)=e^{-I_p(x)}.
	\] For a positive degree-$p$ homogeneous function $J$, layer-cake integration gives
	\[
	\int_{\R^n}e^{-J(x)}\,dx
	=\int_0^\infty e^{-r}|\{J<r\}|\,dr
	=\Gamma(1+n/p)|\{J\le1\}|.
	\]
	Applying the two main theorems and dividing by
$\Gamma(1+n/p)$ gives \eqref{eq:volumes}. For $p\ne2$, the profile
$e^{-r^p}$ is not Gaussian, so Corollary~\ref{cor:radial} gives exactly
the stated equality cases.

At $p=2$, isotropy gives
\[
 H_2(x)=\int_S|P_Ex|^2\,d\nu(E)=|x|^2.
\]
For any section $v$ with barycenter $x$, Lemma~\ref{lem:energy} gives
\[
 \int_S|v(E)|^2\,d\nu(E)\ge|x|^2.
\]
The section $v(E)=P_Ex$ has barycenter $x$ and attains this bound.
Thus $I_2(x)=|x|^2$ as well.
\end{proof}

For $1<p<\infty$, write $p'=p/(p-1)$ for the conjugate exponent. If $K$ is a convex body containing the origin in its interior, write $K^\circ=\{y:\langle x,y\rangle\le1\text{ for all }x\in K\}$ for its polar. Then
\begin{equation}\label{eq:polar-relation}
	K_p^-=(K_{p'}^+)^\circ.
\end{equation}
To verify this, consider the barycenter map on $L^p$ sections. Its quotient norm is $I_p^{1/p}$. Bounded sections give the same infimum: truncate an $L^p$ section and correct its barycenter error $a_R$ by adding $P_Ea_R$; the correction tends to zero in $L^p$. For $\xi\in\R^n$, its dual norm is
\[
 \sup_{\substack{v\in L^p(\nu;\R^n),\ v(E)\in E\ \nu\text{-a.e.}\\
                    \int_S|v(E)|^p\,d\nu(E)\le1}}
 \left\langle\xi,\int_Sv(E)\,d\nu(E)\right\rangle
 =\left(\int_S|P_E\xi|^{p'}\,d\nu(E)\right)^{1/p'}.
\]
This follows either from H\"older's inequality and its equality case or from the dual of a quotient space. Taking polar unit balls gives \eqref{eq:polar-relation}. Thus \eqref{eq:volumes} has the polar form associated with Grassmannian cosine-transform volume inequalities \cite{LXX}; its equality conditions are consequences of the functional classifications.

\subsection{The volume inequalities of Li, Xi, and Zhang}

We now relate Corollary~\ref{cor:volume} to the volume inequalities in
\cite[Theorem 1.2]{LXX}. Let $1\le k<n$, and let $\mu$ be an isotropic
probability measure on $G(n,k)$. Thus $\nu=(n/k)\mu$ satisfies
\eqref{eq:isotropy}. For $1\le p<\infty$, let $Z_{k,p}(\mu)$ be the
origin-symmetric convex body whose support function is
\[
 h_{Z_{k,p}(\mu)}(x)
 =\left(\frac nk\int_{G(n,k)}|P_Ex|^p\,d\mu(E)\right)^{1/p}.
\]
For $p=\infty$, its support function is
\[
 h_{Z_{k,\infty}(\mu)}(x)=\max_{E\in S}|P_Ex|.
\]
Here $S=\operatorname{supp}\mu=\operatorname{supp}\nu$. Write
$Z_{k,p}^*(\mu)$ for the polar body, put
\[
 \alpha_{n,k,\infty}=\omega_k^{n/k},
\]
and use the conventions $1'=\infty$ and $\infty'=1$.

\begin{corollary}\label{cor:lxx-comparison}
For $1\le p\le\infty$,
\begin{equation}\label{eq:lxx-comparison}
 |Z_{k,p}^*(\mu)|\le\alpha_{n,k,p},\qquad
 |Z_{k,p}(\mu)|\ge\alpha_{n,k,p'}.
\end{equation}
For $p=2$, equality holds in both inequalities. For $1<p<\infty$,
$p\ne2$, equality holds in either inequality if and only if $k$ divides
$n$ and
\[
 \mu=\frac{k}{n}\sum_{j=1}^{n/k}\delta_{E_j},
 \qquad \R^n=E_1\oplus\cdots\oplus E_{n/k},
\]
where $E_1,\ldots,E_{n/k}$ are pairwise orthogonal $k$-dimensional
subspaces. The same equality condition holds for the first inequality
when $p=1$ and for the second when $p=\infty$.
\end{corollary}

\begin{proof}
For $1\le p<\infty$, the definitions give
\[
 Z_{k,p}^*(\mu)=K_p^+.
\]
Thus the first inequality follows from Corollary~\ref{cor:volume}, as do
its equality conditions for finite $p$. If $1<p<\infty$, then
\eqref{eq:polar-relation} gives
\[
 Z_{k,p}(\mu)=(K_p^+)^\circ=K_{p'}^-.
\]
The second inequality and its equality conditions for $1<p<\infty$
therefore follow from Corollary~\ref{cor:volume} as well.

As $p\to\infty$, continuity of $(E,x)\mapsto |P_Ex|$ and compactness
of $S$ imply
\[
 h_{Z_{k,p}(\mu)}(x)\longrightarrow
 \max_{E\in S}|P_Ex|=h_{Z_{k,\infty}(\mu)}(x).
\]
The convergence is uniform on the unit sphere: the functions are
equicontinuous there, since their Lipschitz constants are bounded by
$\nu(S)^{1/p}$. Hence the corresponding bodies and their polars
converge in the Hausdorff metric. Similarly, as $p\downarrow1$,
the support functions converge uniformly to $h_{Z_{k,1}(\mu)}$.
The volume is continuous under these convergences. Since
$\alpha_{n,k,p}\to\omega_k^{n/k}$ as $p\to\infty$, the two remaining
endpoint inequalities in \eqref{eq:lxx-comparison} follow.
The equality condition for $|Z_{k,\infty}(\mu)|$ follows from the
finite-exponent equality condition for $K_1^+$ and the relation
$Z_{k,\infty}(\mu)=K_1^-$, which is established below.
\end{proof}

The endpoint relation in the last paragraph requires a separate
argument. For $x\in\R^n$, let
\[
 I_1(x)=\inf_{\substack{v\in L^\infty,\ v(E)\in E\\
 \int v\,d\nu=x}}\int_S|v(E)|\,d\nu(E).
\]
For every such section and every $y\in\R^n$,
\[
 \langle x,y\rangle
 \le \left(\max_{E\in S}|P_Ey|\right)\int_S|v(E)|\,d\nu(E).
\]
Conversely, the closed convex hull of the union of the unit balls
$B_E(1)$ has support function $\max_{E\in S}|P_Ey|$. Every finite convex
combination of points in these balls is a limit of barycenters of
bounded sections whose $L^1(\nu)$ norms do not exceed one: choose small
positive-measure neighborhoods of the relevant points of $S$ and use
the projections of the chosen vectors onto the varying subspaces.
The infimum defining $I_1$ is continuous, so this approximation gives
\[
 K_1^-=Z_{k,\infty}(\mu).
\]

The equality cases at the remaining endpoints, namely the first
inequality in \eqref{eq:lxx-comparison} for $p=\infty$ and the second
for $p=1$, are established in \cite[Theorem 1.2]{LXX} when $k$ divides
$n$. They do not follow merely by taking limits in the finite-exponent
equality conditions.

\appendix
\section{The forward theorem without mass continuity}\label{sec:extension}

The forward inequality also holds under weaker assumptions on the fiber masses. In this appendix, we retain joint continuity of the family and integrability of the positive logarithmic part of the masses, but no longer assume that the masses are continuous. The equality condition is correspondingly stated almost everywhere in the parameter.

\begin{theorem}\label{thm:extended-forward}
	Assume \eqref{eq:isotropy}. Let $f:\cT\to[0,\infty)$ be jointly continuous, with
	\[
	0\le m_E=\int_Ef_E<\infty\quad(E\in S),\qquad
	\int_S\log^+m_E\,d\nu<\infty.
	\]
	Then $\BL(f)\le M:=\exp\int\log m_E\,d\nu\in[0,\infty)$.
	If $M=0$, both sides vanish. If $M>0$, equality holds precisely when
	\eqref{eq:forward-form} holds for almost every parameter and almost every point in that fiber, with measurable positive amplitudes $\theta_E$ and probability densities $h_j$. There is no assumption of mass continuity, a common tail bound, or finite moments.
\end{theorem}
\begin{proof}
	The cutoff argument in Section~\ref{sec:inequalities} already proves the inequality. The logarithmic expressions are well defined: $\log^+m_E$ is integrable by assumption, and the fiber evaluations on a fixed ambient compact set are uniformly bounded. We therefore give the additional argument needed for the equality condition.
	
	Write $u_E(t,y)=f_E*\gamma_{k,t}(y)$ and
	$G_t(x)=\exp\int\log u_E(t,P_Ex)\,d\nu$. The heat family need not be continuous in the parameter. Nevertheless,
	\begin{equation}\label{eq:extended-heat-bounds}
		G_t\ge F_f*\gamma_{n,t},\qquad \int G_t\le M.
	\end{equation}
	For the first relation apply the established forward inequality to
	$f_E(y)\gamma_{k,t}(P_Ex-y)$. Its logarithmic mass positive part is bounded by a constant plus $\log^+m_E$.
	
	For the second, let $f_{E,R}=f_E\chi_R$ and first keep a positive floor:
	\[
	v_{E,R,\varepsilon}(y)
	=(f_{E,R}*\gamma_{k,t})(y)+\varepsilon\gamma_{k,t}(y).
	\]
	This is a jointly continuous family with positive continuous mass
	$m_{E,R}+\varepsilon$, so
	\[
	\int\exp\int\log v_{E,R,\varepsilon}(P_Ex)\,d\nu\,dx
	\le\exp\int\log(m_{E,R}+\varepsilon)\,d\nu.
	\]
	For fixed $\varepsilon>0$, let $R\to\infty$ by monotone convergence, using the finite logarithmic lower bounds supplied by the floor. Next let $\varepsilon\downarrow0$. Subtracting the integrable logarithmic upper bound at $\varepsilon=1$ justifies monotone descent inside the parameter integral. The ambient products are dominated by the integrable product at $\varepsilon=1$. This proves the upper bound in \eqref{eq:extended-heat-bounds}.
	
	Suppose $\int F_f=M>0$. Then $G_t\ge F_f*\gamma_{n,t}>0$ everywhere and the two functions have the same integral. To justify the local differentiations, fix $0<a<b$, $|y|\le K$, and put
	$\ell_E=u_E(a/2,0)$. Since $M>0$, we have $m_E>0$ for almost
	every $E$. On this set $\ell_E>0$, and the Gaussian upper bound gives
	$\ell_E\le C m_E$. Thus
	\[
	 \int_S\log^+\ell_E\,d\nu(E)<\infty,\qquad
	 \int_S\log\ell_E\,d\nu(E)=\log G_{a/2}(0)>-\infty.
	\]
	Consequently $\log\ell_E\in L^1(\nu)$. Comparing Gaussian kernels gives
	\begin{equation}\label{eq:extended-local-bounds}
		c\ell_E\le u_E(t,y)\le C m_E
		\quad(a\le t\le b,\ |y|\le K).
	\end{equation}
	Indeed, $\gamma_{k,t}(y-z)/\gamma_{k,a/2}(z)$ has a positive uniform lower bound in this range.
	
	For $E\in S$ with $m_E>0$, $t>0$, and $y\in E$, let $\pi_{E,t,y}$ be the probability measure on $E$ defined by
 \[
 d\pi_{E,t,y}(z)=
 \frac{f_E(z)e^{-|z-y|^2/(4t)}\,dz}
 {\int_E f_E(w)e^{-|w-y|^2/(4t)}\,dw}.
 \]
The estimates below hold on the full-measure set where $m_E>0$.
	For $r=|z-y|^2/(4t)$ and $s\ge0$, the inequality
	$e^{-r}\mathbf1_{\{r\ge s\}}\le e^{-s/2}e^{-r/2}$, together with
	\eqref{eq:extended-local-bounds}, gives
	\[
	\pi_{E,t,y}\{r\ge s\}
	\le\min\{1,Cm_E\ell_E^{-1}e^{-s/2}\}.
	\]
	Using $\mathbb E_\pi r=\int_0^\infty\pi\{r\ge s\}\,ds$ and
	splitting the integral at $2\log^+(Cm_E/\ell_E)$, we find
	\[
	{\mathbb E}_{\pi}|z-y|^2
	\le C K_E,\qquad
	K_E=1+\log^+m_E+\log^-\ell_E\in L^1(\nu).
	\]
	The explicit identities
	\[
	\nabla\log u_E=\frac{{\mathbb E}_{\pi}z-y}{2t},\quad
	\nabla^2\log u_E=-\frac{I_E}{2t}+\frac{\Cov_\pi(z)}{4t^2},\quad
	\partial_t\log u_E=-\frac{k}{2t}
	+\frac{{\mathbb E}_{\pi}|z-y|^2}{4t^2}
	\]
	therefore give integrable domination of the squared score, second spatial derivatives, and first time derivative on compact sets. Thus $G_t$ is continuous and its almost-everywhere agreement with the convolution becomes the pointwise identity
	\[
	G_t=F_f*\gamma_{n,t}.
	\]
	It is consequently smooth for $t>0$. The defect identity
	\eqref{eq:forward-defect} is now justified without any global integration by parts.
	
	Fix $t>0$. The vanishing of the defect gives, for each $x\in\R^n$,
	\[
	 \nabla_E\log u_E(t,P_Ex)=P_E\nabla\log G_t(x)
	 \quad\text{for }\nu\text{-almost every }E.
	\]
	Intersect the full-measure sets for $x$ in a countable dense subset of $\R^n$. By continuity in $x$, the identity holds for every $x$ on one fixed full-measure set of parameters. Differentiating gives
	\[
	 \nabla^2\log G_t(x)P_E=P_E\nabla^2\log G_t(x)
	 \quad\text{for every }x\in\R^n\text{ and almost every }E.
	\]
	Since $E\mapsto P_E$ is continuous, this commutation relation extends to every $E\in S$.
	
	Apply Lemma~\ref{lem:rigidity} to $\log(G_t/M)$. The rest of the fixed-time recovery in Section~\ref{sec:forward-equality} applies to the probability density $p_0=F_f/M$ and the normalized fibers $f_E/m_E$. It yields the stated almost-everywhere structure. Its sufficiency is the same Fubini and determinant computation as before. The restriction to a full-measure parameter set is retained, since the masses need not be continuous and some exceptional fibers may have zero mass.
	
	More precisely, perform the recovery at one fixed time on the full-measure parameter set just obtained. The global marginals of $p_0$, supplied by Lemma~\ref{lem:deconvolution}, fix the independent factors; the characteristic-function identities of the fibers identify those same factors, so no parameterwise inverse Fourier choice is needed. In the sufficiency direction, a fiber null set has null inverse image under $P_E$ on every bounded ambient cylinder. Fubini therefore transfers the almost-everywhere factorization to the forward product for almost every ambient point. This argument concerns the forward product only; an arbitrary almost-everywhere alteration of a fiber density would not in general be harmless for a supremal envelope.
\end{proof}

\section{A mass-escape obstruction to an unrestricted reverse theorem}\label{sec:escape}

The following example satisfies the assumptions of the forward extension in Appendix~\ref{sec:extension}, but both of its reverse envelopes vanish. It illustrates the role of continuity of the fiber masses in the common admissible class.

Take $n=2$, $k=1$ and let $\nu$ be twice Haar probability measure on $G(2,1)$. Then $\int P_E\,d\nu=I_2$. Fix $e_1=(1,0)$ and $E_0=\R e_1$. For $E\ne E_0$ write $\theta$ for its angle from $E_0$, modulo $\pi$, and set
\[
b_E=\frac{P_Ee_1}{|\sin\theta|},\qquad
f_E(y)=e^{-\pi|y-b_E|^2}.
\]
At $E_0$ define $f_{E_0}\equiv0$. The expression for $b_E$ is independent of the chosen angle representative. It is continuous away from $E_0$, including at the orthogonal line. As $E\to E_0$, $|b_E|\to\infty$, so $f_E$ tends uniformly to zero on bounded fiber sets. Hence $f$ is jointly continuous on the entire tautological bundle.

Every fiber other than $E_0$ has mass one, while $m_{E_0}=0$. Since $\nu$ has no atom at $E_0$,
\[
\exp\int\log m_E\,d\nu=1.
\]
However, $\int|b_E|\,d\nu=\infty$, since $|b_E|$ is comparable to $|\theta|^{-1}$ near zero. For any Bochner-integrable section $v$, finiteness of
$\int|v(E)-b_E|^2\,d\nu$ would imply $v-b\in L^1$ and hence $b\in L^1$, a contradiction. Therefore
\[
\int\log f_E(v(E))\,d\nu
=-\pi\int|v(E)-b_E|^2\,d\nu=-\infty .
\]
Every admissible product is zero, and
\[
\RB(f)\equiv\RBone(f)\equiv0,
\qquad \Mass(f)=1.
\]
Thus joint continuity and logarithmic mass integrability alone do not imply the reverse inequality. This example does not contradict the main theorem: its normalized densities do not extend to a continuous probability family on the compact parameter space. It also explains why boundedness or integrability of the translation section cannot be omitted without replacing it by an appropriate admissibility condition.

\section*{\textbf{Acknowledgments}}
The authors developed the main ideas
and approach of this work and are responsible for the mathematical
results and proofs. ChatGPT 6 was used as an assistive tool in
exploring some proof strategies, checking calculations, and improving
the language and presentation of the manuscript.


\begin{thebibliography}{99}
	
	\bibitem{AS}
	L.~Ambrosio and G.~Savar\'e, \emph{Gradient flows of probability measures},
	lecture notes, March 14, 2007, Section 1.3, pp.~11--12.
	
	\bibitem{Ball}
	K.~Ball, Volume ratios and a reverse isoperimetric inequality,
	\emph{J. London Math. Soc.} (2) \textbf{44} (1991), 351--359.
	
	\bibitem{Barthe}
	F.~Barthe, On a reverse form of the Brascamp--Lieb inequality,
	\emph{Invent. Math.} \textbf{134} (1998), 335--361.
	
	\bibitem{BartheContinuous}
	F.~Barthe, A continuous version of the Brascamp--Lieb inequalities,
	in \emph{Geometric Aspects of Functional Analysis},
	Lecture Notes in Mathematics \textbf{1850}, Springer, 2004, 53--63.
	
	\bibitem{BH}
	F.~Barthe and N.~Huet, On Gaussian Brunn--Minkowski inequalities,
	\emph{Studia Math.} \textbf{191} (2009), 283--304.
	
	\bibitem{BCCT}
	J.~Bennett, A.~Carbery, M.~Christ, and T.~Tao,
	The Brascamp--Lieb inequalities: finiteness, structure and extremals,
	\emph{Geom. Funct. Anal.} \textbf{17} (2008), 1343--1415.
	
	\bibitem{BKX}
	K.~J.~B\"or\"oczky, P.~Kalantzopoulos, and D.~Xi,
	The case of equality in geometric instances of Barthe's reverse Brascamp--Lieb inequality,
	in \emph{Geometric Aspects of Functional Analysis},
	Lecture Notes in Mathematics \textbf{2327}, Springer, 2023, 129--165.
	
	\bibitem{BL}
	H.~J.~Brascamp and E.~H.~Lieb,
	Best constants in Young's inequality, its converse, and its generalization to more than three functions,
	\emph{Adv. Math.} \textbf{20} (1976), 151--173.
	
	\bibitem{BG}
	S.~Brazitikos and A.~Giannopoulos,
	Continuous version of the approximate geometric Brascamp--Lieb inequalities,
	\emph{J. Geom. Anal.} \textbf{32} (2022), article 174.
	
	\bibitem{Brenier}
	Y.~Brenier, Polar factorization and monotone rearrangement of vector-valued functions,
	\emph{Comm. Pure Appl. Math.} \textbf{44} (1991), 365--417.
	
	\bibitem{Caffarelli}
	L.~A.~Caffarelli, Boundary regularity of maps with convex potentials---II,
	\emph{Ann. of Math.} \textbf{144} (1996), 453--496.
	
	\bibitem{CLL}
	E.~A.~Carlen, E.~H.~Lieb, and M.~Loss,
	A sharp analog of Young's inequality on $S^N$ and related entropy inequalities,
	\emph{J. Geom. Anal.} \textbf{14} (2004), 487--520.
	
	\bibitem{CEM}
	D.~Cordero-Erausquin and B.~Maurey,
	Some extensions of the Pr\'ekopa--Leindler inequality using Borell's stochastic approach,
	\emph{Studia Math.} \textbf{238} (2017), 201--233.
	
	\bibitem{CL}
	T.~A.~Courtade and J.~Liu,
	Euclidean forward--reverse Brascamp--Lieb inequalities: finiteness, structure, and extremals,
	\emph{J. Geom. Anal.} \textbf{31} (2021), 3300--3350.
	
	\bibitem{DPF}
	G.~De Philippis and A.~Figalli,
	The Monge--Amp\`ere equation and its link to optimal transportation,
	\emph{Bull. Amer. Math. Soc.} \textbf{51} (2014), 527--580.
	
	\bibitem{FJ}
	A.~Figalli and Y.~Jhaveri,
	Regularity properties of monotone measure-preserving maps,
	\emph{Adv. Nonlinear Stud.} \textbf{23} (2023), article 20220057.
	
	\bibitem{LCH}
	A.-J.~Li, Z.-X.~Cao, and Q.~Huang,
	A matrix form of Grassmannian Ball--Barthe inequality,
	\emph{Linear Multilinear Algebra} \textbf{70} (2022), 6511--6522.
	
	\bibitem{Li2026}
	M.~Li, From the sphere to the Grassmannian: continuous Brascamp--Lieb inequalities via measure approximation and heat flow,
	preprint, 2026.
	arXiv:2610.05890.
	
	\bibitem{LXX}
	A.-J.~Li, D.~Xi, and G.~Zhang,
	Volume inequalities of convex bodies from cosine transforms on Grassmann manifolds,
	\emph{Adv. Math.} \textbf{304} (2017), 494--538.
	
	\bibitem{Lieb}
	E.~H.~Lieb, Gaussian kernels have only Gaussian maximizers,
	\emph{Invent. Math.} \textbf{102} (1990), 179--208.
	
	\bibitem{L1} E. Lutwak, D. Yang, and G. Zhang, \textit{Volume inequalities for subspaces
		of~$L_p$}, J. Differ. Geom. {\bf 68} (2004), 159-184.
	
	
	\bibitem{LYZ}
	E.~Lutwak, D.~Yang, and G.~Zhang,
	Volume inequalities for isotropic measures,
	\emph{Amer. J. Math.} \textbf{129} (2007), 1711--1723.
	
	\bibitem{McCann}
	R.~J.~McCann, Existence and uniqueness of monotone measure-preserving maps,
	\emph{Duke Math. J.} \textbf{80} (1995), 309--323.
	
	\bibitem{TW}
	N.~S.~Trudinger and X.-J.~Wang,
	On the second boundary value problem for Monge--Amp\`ere type equations and optimal transportation,
	\emph{Ann. Sc. Norm. Super. Pisa Cl. Sci.} (5) \textbf{8} (2009), 143--174.
	
	\bibitem{Valdimarsson}
	S.~I.~Valdimarsson, Optimisers for the Brascamp--Lieb inequality,
	\emph{Israel J. Math.} \textbf{168} (2008), 253--274.
	
\end{thebibliography}
\end{document}